\documentclass[11pt]{article}

\usepackage[english]{babel}
\usepackage[utf8]{inputenc}
\usepackage[T1]{fontenc}

\usepackage{amsthm}
\usepackage{amssymb}
\usepackage{amsfonts}
\usepackage{amsmath}
\usepackage{mathtools}

\usepackage{tikz}
\usepackage{float}

\usepackage[top=20mm, bottom=20mm, left=20mm, right=20mm, a4paper]{geometry}
\usepackage[none]{hyphenat}

\usepackage[colorlinks=true, citecolor=red, linkcolor=blue, urlcolor=blue]{hyperref}

\allowdisplaybreaks

\newtheorem{theorem}{Theorem}

\newtheorem{lemma}{Lemma}
\newtheorem{proposition}{Proposition}
\newtheorem{remark}{Remark}

\theoremstyle{definition}
\newtheorem{definition}{Definition}

\numberwithin{equation}{section}

\makeatletter
\renewcommand{\maketitle}{%
  \begin{center}
    {\Large\bfseries \@title \par}
    \vskip 1.2em
    {\normalsize \@author \par}
    \vskip 0.8em
  \end{center}
}
\makeatother

\makeatletter
\renewenvironment{proof}[1][Proof]{\par
    \pushQED{\qed}%
    \normalfont \topsep6\p@\@plus6\p@\relax
    \trivlist
    \item\relax
    {\textbf{#1\@addpunct{.}}}\hspace\labelsep\ignorespaces
}{%
    \popQED\endtrivlist\@endpefalse
}
\makeatother

\let\OLDthebibliography\thebibliography
\renewcommand{\thebibliography}[1]{%
    \OLDthebibliography{#1}%
    \small
    \setlength{\parskip}{0pt}%
    \setlength{\itemsep}{0pt plus 0.3ex}%
}

\title{Nonzero solutions for Dirichlet problems\linebreak with the competing
\texorpdfstring{$\boldsymbol{(p,q)}$}{(p,q)}-Laplacian}

\author{Micha{\l} Be{\l}dzi\'nski}

\date{}

\begin{document}

\sloppy

\maketitle

\begin{center}
\small
Institute of Mathematics\\
Lodz University of Technology\\
al. Politechniki 8\\
93-590 Lodz, Poland\\
ORCID: \href{https://orcid.org/0000-0002-9653-612X}{0000-0002-9653-612X}\\
\texttt{michal.beldzinski@p.lodz.pl}
\end{center}

\begin{abstract}
    The following equation is considered:
    \begin{equation*}
        \begin{cases}
            -\Delta_p u(\boldsymbol{x}) + \Delta_q u(\boldsymbol{x}) = f\big(\boldsymbol{x}, u(\boldsymbol{x})\big)
            & \text{for } \boldsymbol{x} \in \Omega, \\
            u(\boldsymbol{x}) = 0
            & \text{for } \boldsymbol{x}\in \partial \Omega.
        \end{cases}
    \end{equation*}
    We start by showing that the energy functional associated with the considered problem no longer has the properties typical of its non-competing counterpart, which prevents a direct application of standard variational and monotonicity methods. Using Ekeland's variational principle combined with properties imposed on the nonlinear term~$f$, we prove the existence of a nonzero generalized solution to the problem. It is also shown that in the one-dimensional case $\Omega = (-1,1)$, assuming that $f$ satisfies the symmetry conditions
    \begin{equation*}
        f(x,u) = f(-x,u) = -f(x,-u),
    \end{equation*}
    together with the sign condition $f(x,u) \geq 0$ for $u \geq 0$ and an appropriate growth bound on $f$, a convex relaxation of the action functional yields the existence of a positive weak solution. Furthermore, we prove that in this setting, and under a slightly stronger sign condition
    \begin{equation*}
        f(x,u) > 0\text{ for }u > 0,
    \end{equation*}
    all global minimizers of the action functional do not vanish in $(-1,1)$ but fail to be of class~$C^1$.
\end{abstract}

\section{Introduction}

The competing $(p,q)$-Laplacian $-\Delta_p + \Delta_q$ originally appeared in \cite{Liu2020, Motreanu2020}, where it serves as a natural extension of the weighted $(p,q)$-Laplacian operator $-\Delta_p + \mu\Delta_q$, with $\mu$ being a constant weight. In the case $\mu > 0$, the operator $-\Delta_p + \mu\Delta_q$ loses its monotonicity and the crucial characteristics needed for using standard techniques to show the existence of solutions. One method to address the aforementioned difficulty is to introduce a new, more comprehensive notion of solution called a \emph{generalized solution}. In the existing literature on competing operators, the analysis is formulated primarily in terms of generalized solutions, while their relation to classical weak solutions has remained largely unexamined.
\par
The competing $(p,q)$-Laplacian has driven a wide range of mathematical investigations, including differential inclusions \cite{Cen2025}, equations with unbounded weights \cite{Diblik2025, Galewski2024}, fractional operators \cite{Gambera2024}, and convolution effects \cite{Motreanu2021, Razani2022}. Authors have also considered Heisenberg group settings \cite{Razani2025, Razani2023}, gradient-dependent nonlinearities \cite{Razani2023a}, logarithmic-type double phase problems \cite{Zhong2025}, as well as recent extensions to variable exponent or multi-phase frameworks \cite{Vetro2026, Zhang2026}. While these works show that generalized solutions provide a flexible framework for competing operators, they leave open the question whether classical weak solutions can be obtained under additional structural assumptions. A comprehensive survey book \cite{Galewski2026} encapsulates these diverse advancements in nonlinear analysis.
\par
The present paper is devoted to problems involving the competing $(p,q)$-Laplacian of the form
\begin{equation}
    \label{eq:main_problem}
    \begin{cases}
        -\Delta_p u(\boldsymbol{x}) + \Delta_q u(\boldsymbol{x}) = f\big(\boldsymbol{x}, u(\boldsymbol{x})\big) & \text{for } \boldsymbol{x} \in \Omega, \\
        u(\boldsymbol{x}) = 0 & \text{for }\boldsymbol{x}\in \partial \Omega,
    \end{cases}
\end{equation}
where ${\Omega\subset \mathbb{R}^d}$ is open, bounded and non-empty, and ${1<q<p<\infty}$. The nonlinear term ${f\colon \Omega \times \mathbb{R}\longrightarrow \mathbb{R}}$ is a Carath{\'e}odory function satisfying suitable growth and symmetry conditions specified later in the paper. Because the associated energy functional lacks weak lower semicontinuity and the Palais--Smale condition, and its derivative lacks pseudomonotonicity, standard tools of nonlinear analysis cannot be applied directly. We will show that under a standard assumption concerning the growth of~$f$ and under the assumption that the antiderivative of $f(\boldsymbol{x},\cdot)$ is positive on a right neighbourhood of $0$, using Ekeland's variational principle allows us to prove that problem \eqref{eq:main_problem} has a nonzero generalized solution.
\par
In the one-dimensional case, however, we are able to establish the existence of nonzero weak solutions. More precisely, we consider the boundary value problem
\begin{equation}
    \label{eq:main_problem-1dim}
    \left\{
    \begin{array}{cl}
        -\left(|u'(x)|^{p-2}u'(x) - |u'(x)|^{q-2}u'(x)\right)' = f\big(x, u(x)\big) & \text{for } x \in (-1,1), \\
        u(-1) = u(1) = 0. &
    \end{array}
    \right.
\end{equation}
Using a convex relaxation of the action functional together with symmetry and sign conditions on $f$, we show that problem \eqref{eq:main_problem-1dim} admits a positive weak solution, even though the associated energy functional has no global minimizers of class $C^1$. This relaxation approach is related to the method used in \cite{Kromer} for autonomous radially symmetric variational problems on a ball, where the nonlinearity depends solely on $u$. In our one-dimensional setting, an explicit modification of the derivative allows us to treat nonlinearities $f$ that depend explicitly on the spatial variable $x$.
\par
The paper is organized as follows. Section \ref{sec:preliminaries} reviews the basics of nonlinear functional analysis used for proving the existence of nonzero solutions for potential nonlinear equations and provides an auxiliary construction of arbitrarily small functions with unit gradient almost everywhere. In Section \ref{sec:energy_properties}, it will be shown that none of the properties recalled above are satisfied by the functional associated with the main equation. This serves as a motivation for considering the so-called generalized solution. The definition of this solution is presented in Section \ref{sec:generalized_solution}, along with the corresponding existence result. In Section \ref{sec:weak_solution}, it will be demonstrated that in $\Omega = (-1,1)$, by strengthening the assumptions on the nonlinearity $f$, a weak solution exists for any $1 < q < p < \infty$. It will also be shown that this solution is not of class $C^1$. The obtained result will be illustrated with a typical example and compared with its multidimensional autonomous counterpart on a ball via \cite{Kromer}.

\section{Preliminaries}
\label{sec:preliminaries}

We begin by recalling selected concepts from nonlinear functional analysis that are used in standard methods for proving the existence of solutions. Let $(X, \Vert{}{\cdot}\Vert{})$ be a reflexive Banach space. A functional $\mathcal{I} \colon X \longrightarrow \mathbb{R}$ is called \emph{coercive} if
\begin{equation*}
    \lim_{\|u\| \to \infty} \mathcal{I}(u) = \infty.
\end{equation*}
Coercivity ensures the existence of a bounded minimizing sequence. In infinite-dimensional Banach spaces, however, it typically has to be supported by an additional assumption in order to obtain the existence of a minimizer. One possibility is closely related to the weak topology. Recall that we say that $\mathcal{I}$ is \emph{w.l.s.c.} (weakly lower semicontinuous) if $\{ v \in X : \mathcal{I}(v) \leq \alpha \}$ is weakly closed for each $\alpha \in \mathbb{R}$, and \emph{s.w.l.s.c.} (sequentially weakly lower semicontinuous) if 
\begin{equation*}
    v_n \rightharpoonup v\text{ in }X \implies \liminf_{n \to \infty} \mathcal{I}(v_n) \geq \mathcal{I}(v).
\end{equation*}
Using the weak compactness of the closed unit ball in reflexive Banach spaces, one can prove the following result.

\begin{proposition}[{\cite[Theorem 2.22]{Galewski2021}}]
    \label{prop:DirectMethod}
    Let $\mathcal{I}$ be a coercive $C^1$ functional defined on a reflexive Banach space $X$. If $\mathcal{I}$ is s.w.l.s.c., then there exists $u \in X$ satisfying
    \begin{equation*}
        \mathcal{I}(u) = \min_{v \in X} \mathcal{I}(v).
    \end{equation*}
    Such $u$ is hence a critical point of $\mathcal{I}$.
\end{proposition}

To shift the burden of ensuring compactness from the space to the functional itself, one can employ Ekeland's variational principle, which we state in the version for $C^1$ functionals following \cite[Corollary~3.3]{Jabri2003}.

\begin{theorem}[the Ekeland Variational Principle]
    \label{thm:Ekeland}
    Let $X$ be a Banach space. If $\mathcal{I} \colon X \longrightarrow \mathbb{R}$ is a $C^1$ functional bounded from below, then there exists a sequence $(u_n) \subset X$ satisfying
    \begin{equation*}
        \lim_{n\to \infty}\mathcal{I}(u_n) = \inf_{v \in X} \mathcal{I}(v) \quad \text{and}\quad \lim_{n\to \infty} \mathcal{I}'(u_n) = 0.
    \end{equation*}
\end{theorem}

The condition that is the most commonly used to complement Ekeland's variational principle is \emph{(PS)} (\emph{the Palais--Smale condition}, see \cite{Mawhin2010}): every sequence $(x_n) \subset X$ satisfying $\sup_{n \in \mathbb{N}} |\mathcal{I}(x_n)| < \infty$ and $\mathcal{I}'(x_n) \to 0$ has a convergent subsequence. One can easily show that a $C^1$ functional bounded from below that satisfies the Palais--Smale condition admits a global minimizer. For a more detailed discussion we refer the reader to \cite[Chapter 3]{Jabri2003}.
\par
Another approach to study the solvability of
\begin{equation*}
    \mathcal{I}'(u) = 0
\end{equation*}
is connected with monotonicity methods. We now recall a standard and relatively non-restrictive assumption imposed on the operator within this framework. An operator ${A \colon X \longrightarrow X^*}$ (in particular, ${A = \mathcal{I}'}$) is called \emph{pseudomonotone} if
\begin{equation*}
    \left.
    \begin{array}{c}
        u_n \rightharpoonup u \\
        \limsup_{n \to \infty} \langle A(u_n), u_n - u \rangle \leq 0
    \end{array}
    \right\}\implies \liminf_{n \to \infty} \langle A(u_n), u_n - v \rangle \geq \langle A(u), u - v \rangle\text{ for every }v \in X.
\end{equation*}
This condition is met by most nonlinear operators such as the $p$-Laplacian (with ${1 < p < \infty}$) or its counterparts like double phase operators; see, for instance, \cite[Section 3.4]{Migorski2013} and \cite{Papageorgiou2019} for more examples and a detailed discussion on pseudomonotonicity. Let us also note that the above assumption is satisfied by any monotone operator, provided we impose a minimal continuity condition, for instance hemicontinuity or radial continuity; see \cite{Galewski2021}.
\par
Another branch of assumptions imposed on $A$ consists of various versions of the (S) condition. The~least restrictive one is recalled here following \cite{Galewski2021}. It is said that $A$ satisfies \emph{the condition (S)$_2$} if
\begin{equation*}
    \left.
    \begin{array}{c}
        u_n \rightharpoonup u\\
        A(u_n) \to A(u)
    \end{array}
    \right\} \implies u_n \to u.
\end{equation*}
Following \cite{Leoni2017}, the first-order Sobolev space on an open, bounded and non-empty set $\Omega \subset \mathbb{R}^d$ is denoted by ${W^{1,p}(\Omega)}$ for $p\in [1,\infty]$. Recall that
\begin{equation*}
    W^{1,p}_0(\Omega) \coloneq \overline{C_c^\infty(\Omega)}^{\,W^{1,p}(\Omega)}.
\end{equation*}
Moreover, for $1 \leq p < \infty$, there exists a constant ${C_p > 0}$ such that the following \emph{Poincar{\'e} inequality} holds:
\begin{equation}
    \label{eq:Poincare}
    \|u\|_{W^{1,p}_0}^p = \int_\Omega |\nabla u(\boldsymbol{x})|^p \, d\boldsymbol{x} \geq C_p \int_\Omega |u(\boldsymbol{x})|^p \, d\boldsymbol{x} = C_p \|u\|_{L^p}^p \quad \text{for every } u \in W^{1,p}_0(\Omega).
\end{equation}
Throughout the paper, we take $C_p$ to be the optimal (i.e.\ the largest possible) constant in \eqref{eq:Poincare}.\linebreak For $1 < p < \infty$, we denote the normed dual of $W^{1,p}_0(\Omega)$ by $W^{-1,p'}(\Omega)$, where $p' = \frac{p}{p-1}$ stands for the H\"older-conjugate exponent.
\par
By a cube centred at $\boldsymbol{y} = (y_1,\ldots, y_d)\in \mathbb{R}^d$ with half-side length $r > 0$ we mean any set of the form $Q = \boldsymbol{y} + [-r,r]^d$. For every such set we define a function $\delta_Q \colon \mathbb{R}^d \longrightarrow \mathbb{R}$ by the formula
\begin{equation*}
    \delta_Q(\boldsymbol{x}) \coloneq \operatorname{dist}(\boldsymbol{x}, \mathbb{R}^d \setminus Q)
    \qquad \text{for every }\boldsymbol{x}\in \mathbb{R}^d.
\end{equation*}
We will identify $\delta_Q$ with its restriction $\delta_Q|_\Omega$. Under that convention, since $\delta_Q$ is Lipschitz continuous on $\mathbb{R}^d$ with Lipschitz constant $1$, Rademacher's theorem (see \cite[Theorem 3.2]{EvansGariepy2015}) implies that $\delta_Q \in W^{1,\infty}(\Omega)$. Moreover:
\begin{itemize}
    \item $0 \leq \delta_Q(\boldsymbol{x})\leq r$ for every $\boldsymbol{x}\in \mathbb{R}^d$;
    \item $\left\vert{}\nabla \delta_Q(\boldsymbol{x})\right\vert{} = 1$ for a.e. $\boldsymbol{x}\in Q$;
    \item $\nabla \delta_Q(\boldsymbol{x}) = \boldsymbol{0}$ for a.e. $\boldsymbol{x} \in \mathbb{R}^d \setminus Q$
\end{itemize}
since for every $\boldsymbol{x} = (x_1,\ldots, x_d) \in \mathbb{R}^d$ we have
\begin{equation}
    \label{eq:distance_gradient}
     \delta_Q(\boldsymbol{x}) = 
     \begin{cases}
         \displaystyle \min_{1\leq i \leq d} \big(r - |y_i - x_i|\big) & \text{if }\boldsymbol{x}\in Q,\\
         0 & \text{otherwise}.
     \end{cases}
\end{equation}

\begin{figure}[!htbp]
    \centering
    \begin{tikzpicture}[
        x={(-0.45cm,-0.32cm)},
        y={(1cm,0cm)},
        z={(0cm,1cm)},
        >=stealth,
        scale=1.05
    ]
        \coordinate (A) at (-2,-2,0);
        \coordinate (B) at ( 2,-2,0);
        \coordinate (C) at ( 2, 2,0);
        \coordinate (D) at (-2, 2,0);
        \coordinate (Y) at ( 0, 0,0);
        \coordinate (P) at ( 0, 0,2);

        \draw[->, thin, gray!80!black] (-3.2, 0, 0) -- (3.4, 0, 0) node[below left, black] {$x_1$};
        \draw[->, thin, gray!80!black] (0, -3.0, 0) -- (0, 3.4, 0) node[right, black] {$x_2$};
        \draw[->, thin, gray!80!black] (0, 0, 0) -- (0, 0, 2.9) node[above, black] {$\delta_Q(\boldsymbol{x})$};

        \fill[gray!15, opacity=0.85] (B) -- (C) -- (P) -- cycle;
        \fill[gray!28, opacity=0.85] (C) -- (D) -- (P) -- cycle;

        \draw[dashed, gray!70!black] (D) -- (A) -- (B);
        \draw[dashed, gray!70!black] (A) -- (P);

        \draw[dotted, gray!80!black] (A) -- (C);
        \draw[dotted, gray!80!black] (B) -- (D);

        \draw[dashed, thick, black] (Y) -- (P) node[midway, right=2pt, font=\footnotesize] {$r$};
        \node[circle, fill=black, inner sep=1.1pt, label={[below left=-1pt, font=\footnotesize]$\boldsymbol{y}$}] at (Y) {};
        \node[circle, fill=black, inner sep=1.2pt] at (P) {};

        \draw[very thick] (B) -- (C) -- (D);
        \draw[very thick] (B) -- (P);
        \draw[very thick] (C) -- (P);
        \draw[very thick] (D) -- (P);

        \node[below right, font=\footnotesize] at (1.1, 2.0, 0) {$Q = \boldsymbol{y} + [-r,r]^2$};
        \draw[<->, thin, gray!80!black] (2.25, 0, 0) -- (2.25, 2, 0) node[midway, below right=-1pt, font=\footnotesize, black] {$r$};
    \end{tikzpicture}
    \caption{Graph of the function $\delta_Q$ over a cube $Q = \boldsymbol{y} + [-r,r]^2$.}
    \label{fig:delta_Q}
\end{figure}
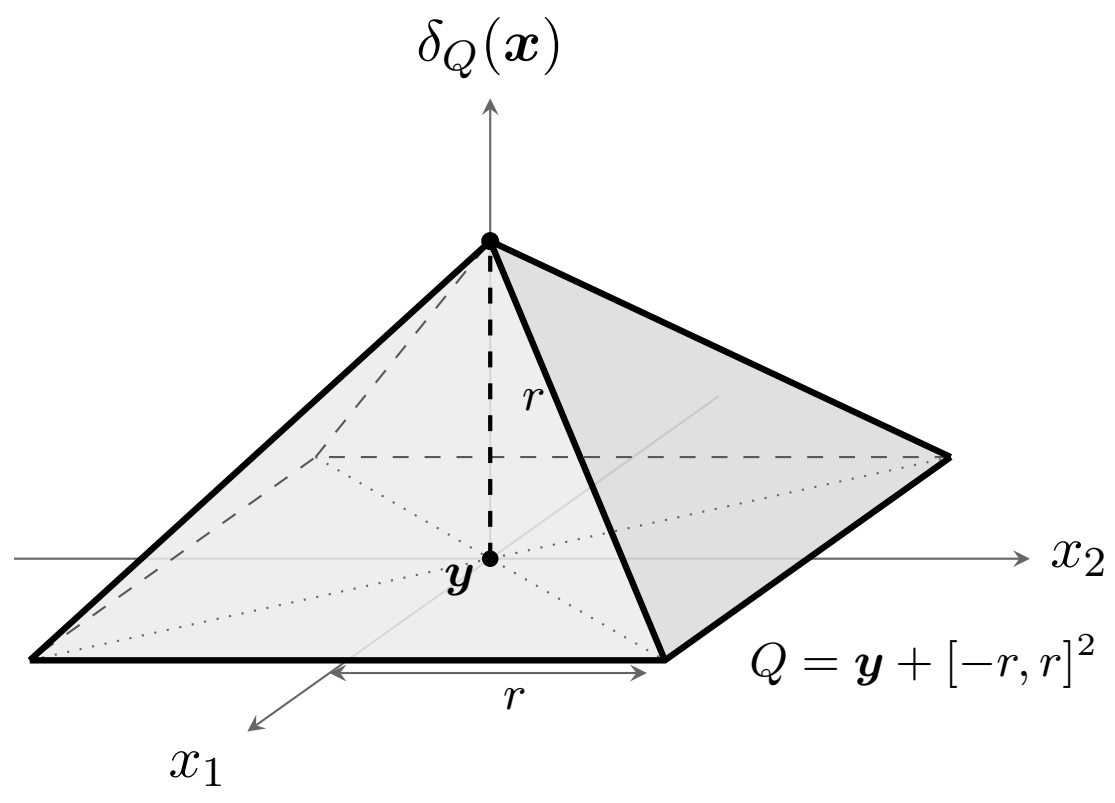

Notice also that if $Q \subset \Omega$, then $\operatorname{supp}(\delta_Q) = Q$ is a compact subset of $\Omega$, which implies $\delta_Q \in W^{1,p}_0(\Omega)$ for every $1 \leq p < \infty$. Based on the functions $\delta_Q$, we can prove the following lemma.

\begin{lemma}
    \label{lem:unit_gradient_sequence}
    Take a non-empty, open and bounded set $\Omega \subset \mathbb{R}^d$ and let $1 < p < \infty$. There exists a~sequence $(u_n) \subset W^{1,p}_0(\Omega)$ satisfying:
    \begin{itemize}
        \item $u_n(\boldsymbol{x}) > 0$ for all $n\in \mathbb{N}_+$ and a.e. $\boldsymbol{x} \in \Omega$,
        \item $|\nabla u_n(\boldsymbol{x})| = 1$ for all $n\in \mathbb{N}_+$ and a.e. $\boldsymbol{x}\in \Omega$,
        \item $\|u_n\|_\infty \leq \frac{1}{2n}$ for every $n\in \mathbb{N}_+$,
        \item $u_n \rightharpoonup 0$ in $W^{1,p}_0(\Omega)$.
    \end{itemize}
\end{lemma}
\begin{proof}
    Since $\Omega \subset \mathbb{R}^d$ is open, bounded and non-empty, its complement $\mathbb{R}^d \setminus \Omega$ is closed and non-empty. By the Whitney decomposition theorem (see \cite[Appendix J, Proposition J.1]{Grafakos2014}), there exists a countable family $\mathcal{Q} = \{Q_j\}_{j=1}^\infty$ of closed cubes $Q_j = \boldsymbol{y}_j + [-r_j, r_j]^d$ with centres $\boldsymbol{y}_j \in \Omega$ and half-side lengths $r_j > 0$ such that:
    \begin{enumerate}
        \item[(i)] $\Omega = \bigcup_{j=1}^\infty Q_j$;
        \item[(ii)] the cubes $Q_j$ have pairwise disjoint interiors, that is, $\operatorname{int}(Q_j) \cap \operatorname{int}(Q_l) = \emptyset$ whenever $j \neq l$;
        \item[(iii)] for every $j \in \mathbb{N}_+$, the diameter $\operatorname{diam}(Q_j) = 2\sqrt{d}\,r_j$ satisfies
        \begin{equation*}
            \operatorname{diam}(Q_j) \leq \operatorname{dist}(Q_j, \partial \Omega) \leq 4\operatorname{diam}(Q_j).
        \end{equation*}
    \end{enumerate}
    
    \begin{figure}[!htbp]
        \centering
        \begin{tikzpicture}[scale=0.82]
            \fill[gray!6] (0,0) ellipse (4.6cm and 3.3cm);
            \draw[very thick] (0,0) ellipse (4.6cm and 3.3cm);
            \node[above right, font=\small] at (3.25, 2.35) {$\partial\Omega$};
            \node[font=\small, fill=white, inner sep=1.5pt] at (-1.0, 1.0) {$Q_j$};

            \foreach \x/\y/\s in {
                0/0/2, 2/0/1, 3/0/0.5, 3.5/0/0.5, 3/0.5/0.5, 3.5/0.5/0.5, 2/1/0.5, 2.5/1/0.5,
                2/1.5/0.5, 2.5/1.5/0.5, 3/1/0.5, 3.5/1/0.25, 3.75/1/0.25, 3.5/1.25/0.25, 3.75/1.25/0.125, 3.875/1.25/0.125,
                3.75/1.375/0.125, 3.875/1.375/0.125, 3/1.5/0.25, 3.25/1.5/0.25, 3/1.75/0.25, 3.25/1.75/0.25, 3.5/1.5/0.25, 3.75/1.5/0.125,
                3.5/1.75/0.125, 0/2/0.5, 0.5/2/0.5, 0/2.5/0.5, 0.5/2.5/0.5, 1/2/0.5, 1.5/2/0.5, 1/2.5/0.25,
                1.25/2.5/0.25, 1/2.75/0.25, 1.25/2.75/0.25, 1.5/2.5/0.25, 1.75/2.5/0.25, 1.5/2.75/0.125, 1.625/2.75/0.125, 1.5/2.875/0.125,
                1.75/2.75/0.125, 1.875/2.75/0.125, 0/3/0.125, 0.125/3/0.125, 0.25/3/0.125, 0.375/3/0.125, 0.5/3/0.125, 0.625/3/0.125,
                0.75/3/0.125, 0.875/3/0.125, 1/3/0.125, 2/2/0.5, 2.5/2/0.25, 2.75/2/0.25, 2.5/2.25/0.25, 2.75/2.25/0.125,
                2.875/2.25/0.125, 2.75/2.375/0.125, 2/2.5/0.25, 2.25/2.5/0.125, 2.375/2.5/0.125, 2.25/2.625/0.125, 2.5/2.5/0.125, 3/2/0.125,
                3.125/2/0.125, 3/2.125/0.125, 3.125/2.125/0.125, 3.25/2/0.125, 4/0/0.25, 4.25/0/0.125, 4.375/0/0.125, 4.25/0.125/0.125,
                4.375/0.125/0.125, 4/0.25/0.25, 4.25/0.25/0.125, 4.375/0.25/0.125, 4.25/0.375/0.125, 4/0.5/0.25, 4.25/0.5/0.125, 4.25/0.625/0.125,
                4/0.75/0.25, 4.25/0.75/0.125, 4/1/0.125, 4.125/1/0.125, 4/1.125/0.125
            } {
                \draw[thin, fill=white, fill opacity=0.6] ( \x,  \y) rectangle ( \x+\s,  \y+\s);
                \draw[thin, fill=white, fill opacity=0.6] (-\x,  \y) rectangle (-\x-\s,  \y+\s);
                \draw[thin, fill=white, fill opacity=0.6] ( \x, -\y) rectangle ( \x+\s, -\y-\s);
                \draw[thin, fill=white, fill opacity=0.6] (-\x, -\y) rectangle (-\x-\s, -\y-\s);
            }
            \fill[gray!25, draw=black, thick] (-2, 0) rectangle (0, 2);
            \node[font=\small] at (-1.0, 1.0) {$Q_j$};
        \end{tikzpicture}
        \caption{Illustration of the Whitney decomposition $\mathcal{Q} = \{Q_j\}_{j=1}^\infty$ of a bounded open set $\Omega \subset \mathbb{R}^2$.}
        \label{fig:Whitney_decomposition}
    \end{figure}
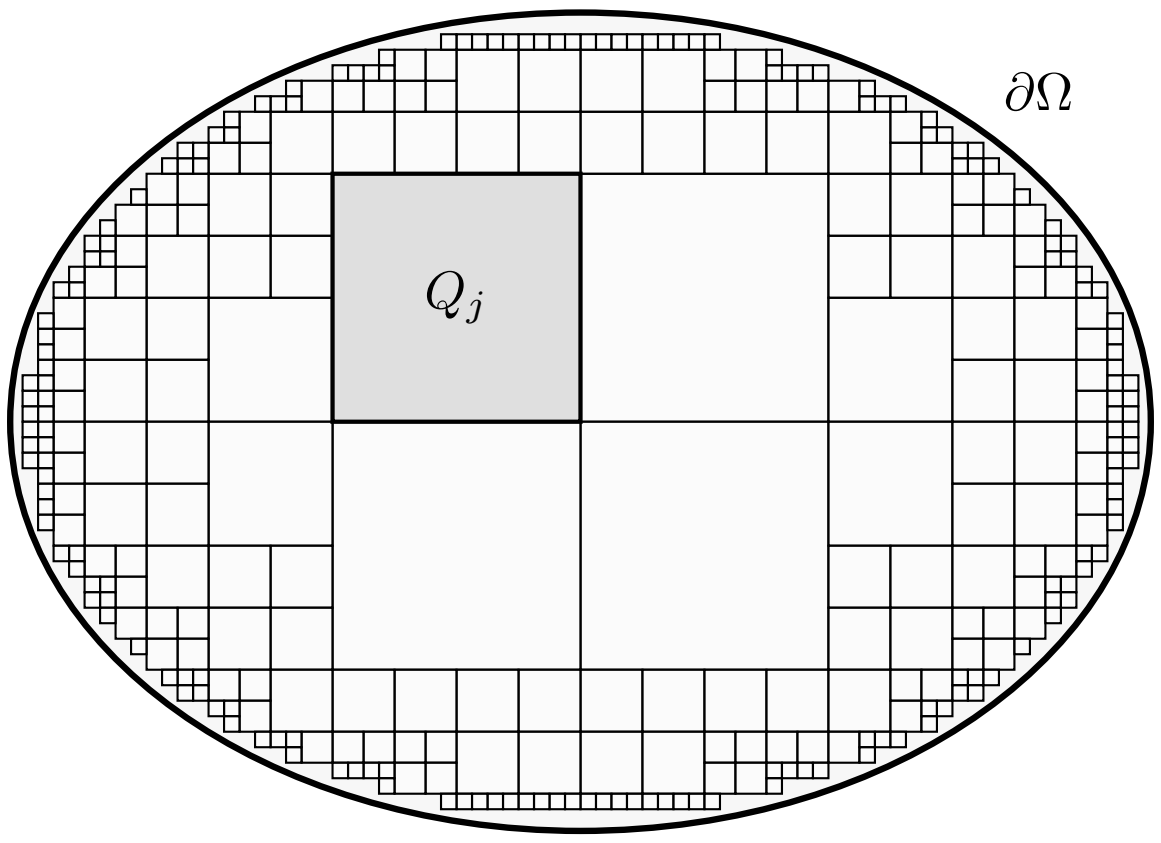
    
    Since $\Omega$ is bounded, we have
    \begin{equation*}
        R_{\mathcal{Q}} \coloneq \max_{j \in \mathbb{N}_+} r_j \leq \tfrac{\operatorname{diam}(\Omega)}{2\sqrt{d}} < \infty.
    \end{equation*}
    Notice also that since the boundary $\partial Q_j$ of each cube has $d$-dimensional Lebesgue measure zero, the~union $\bigcup_{j=1}^\infty \partial Q_j$ is a null set, and hence $\left|\Omega \setminus \bigcup_{j=1}^\infty \operatorname{int}(Q_j)\right| = 0$. Therefore, from condition (ii) it~follows that
    \begin{equation*}
        \sum_{j=1}^\infty |Q_j| = \sum_{j=1}^\infty \left|\operatorname{int}(Q_j)\right| = \left|\bigcup_{j=1}^\infty \operatorname{int}(Q_j)\right| = |\Omega| < \infty.
    \end{equation*}
    For each $j \in \mathbb{N}_+$, consider the function $\delta_{Q_j}$ defined in \eqref{eq:distance_gradient}. Since $Q_j \subset \Omega$ is compact, $\delta_{Q_j}$ is Lipschitz continuous on $\mathbb{R}^d$, and $\operatorname{supp}(\delta_{Q_j}) = Q_j \subset \Omega$, we have $\delta_{Q_j} \in W^{1,p}_0(\Omega)$ (see \cite[Theorem 11.35]{Leoni2017}). Consider the pointwise series
    \begin{equation*}
        w_{\mathcal{Q}}(\boldsymbol{x}) \coloneq \sum_{j=1}^\infty \delta_{Q_j}(\boldsymbol{x}) \qquad \text{for } \boldsymbol{x} \in \Omega.
    \end{equation*}
    To see that $w_{\mathcal{Q}}$ is well defined it is enough to observe that for any $n, m \in \mathbb{N}_+$, $n < m$, we have
    \begin{equation*}
        \left\|\sum_{j = n}^m \delta_{Q_j}\right\|_{W^{1,p}_0} = \left( \int_\Omega \left| \sum_{j = n}^m \nabla \delta_{Q_j}(\boldsymbol{x}) \right|^p d\boldsymbol{x} \right)^{\frac{1}{p}} = \left(\sum_{j = n}^m |Q_j|\right)^{\frac{1}{p}} \leq \left(\sum_{j = n}^\infty |Q_j|\right)^{\frac{1}{p}}
    \end{equation*}
    since the interiors $\left\{\operatorname{int}Q_j\right\}_{j = 1}^\infty$ are pairwise disjoint sets. Moreover ${\left\|\delta_{Q_j}\right\|_\infty = r_j}$ and hence 
    \begin{equation*}
        \left\|w_{\mathcal{Q}}\right\|_\infty = R_{\mathcal{Q}}.
    \end{equation*}
    Finally, $w_{\mathcal{Q}}(\boldsymbol{x}) > 0$ if and only if $\boldsymbol{x}\in \bigcup_{j = 1}^\infty \operatorname{int}Q_j$, which shows that $w_{\mathcal{Q}}(\boldsymbol{x}) > 0$ for a.e. $\boldsymbol{x}\in \Omega$. 
    
    Now, each cube $Q = \boldsymbol{y} + [-r,r]^d$ can be divided into a family of $3^d$ closed subcubes of equal size. If we denote by $\{\boldsymbol{e}_1,\ldots, \boldsymbol{e}_d\}$ the canonical basis in $\mathbb{R}^d$, the explicit formula for such a division reads
    \begin{equation*}
        D(Q) = \bigcup_{\boldsymbol{k} \in \{-1,0,1\}^d}\left\{\boldsymbol{y} + \tfrac{2r}{3}\sum_{j = 1}^d k_j \boldsymbol{e}_j + \left[-\tfrac{r}{3}, \tfrac{r}{3}\right]^d \right\}.
    \end{equation*}
    
    \begin{figure}[!htbp]
        \centering
        \begin{tikzpicture}[
            x=2.5cm, y=2.5cm,
            >=stealth,
            cell/.style={
                below=4pt,
                align=center,
                fill=gray!6,
                inner sep=1.5pt,
                rounded corners=1pt
            },
            centercell/.style={
                below=4pt,
                align=center,
                fill=gray!20,
                inner sep=1.5pt,
                rounded corners=1pt
            },
            dot/.style={circle, fill=black, inner sep=1.2pt}
        ]
            {\footnotesize
            \fill[gray!6] (0,0) rectangle (3,3);
            \fill[gray!20] (1,1) rectangle (2,2);

            \draw[dotted, gray!80!black] (0.5, 0.5) -- (2.5, 0.5);
            \draw[dotted, gray!80!black] (0.5, 1.5) -- (2.5, 1.5);
            \draw[dotted, gray!80!black] (0.5, 2.5) -- (2.5, 2.5);
            \draw[dotted, gray!80!black] (0.5, 0.5) -- (0.5, 2.5);
            \draw[dotted, gray!80!black] (1.5, 0.5) -- (1.5, 2.5);
            \draw[dotted, gray!80!black] (2.5, 0.5) -- (2.5, 2.5);

            \draw[thick, gray!80!black] (1,0) -- (1,3) (2,0) -- (2,3);
            \draw[thick, gray!80!black] (0,1) -- (3,1) (0,2) -- (3,2);
            \draw[very thick, black] (0,0) rectangle (3,3);

            \node[dot] at (0.5, 2.5) {};
            \node[cell] at (0.5, 2.5) {$\boldsymbol{y} - \tfrac{2r}{3}\boldsymbol{e}_1 + \tfrac{2r}{3}\boldsymbol{e}_2$ \\[3pt] $+ \left[-\tfrac{r}{3}, \tfrac{r}{3}\right]^d$};

            \node[dot] at (1.5, 2.5) {};
            \node[cell] at (1.5, 2.5) {$\boldsymbol{y} + 0\boldsymbol{e}_1 + \tfrac{2r}{3}\boldsymbol{e}_2$ \\[3pt] $+ \left[-\tfrac{r}{3}, \tfrac{r}{3}\right]^d$};

            \node[dot] at (2.5, 2.5) {};
            \node[cell] at (2.5, 2.5) {$\boldsymbol{y} + \tfrac{2r}{3}\boldsymbol{e}_1 + \tfrac{2r}{3}\boldsymbol{e}_2$ \\[3pt] $+ \left[-\tfrac{r}{3}, \tfrac{r}{3}\right]^d$};

            \node[dot] at (0.5, 1.5) {};
            \node[cell] at (0.5, 1.5) {$\boldsymbol{y} - \tfrac{2r}{3}\boldsymbol{e}_1 + 0\boldsymbol{e}_2$ \\[3pt] $+ \left[-\tfrac{r}{3}, \tfrac{r}{3}\right]^d$};

            \node[dot, label={[above right=0pt, font=\footnotesize]$\boldsymbol{y}$}] at (1.5, 1.5) {};
            \node[centercell] at (1.5, 1.5) {$\boldsymbol{y} + 0\boldsymbol{e}_1 + 0\boldsymbol{e}_2$ \\[3pt] $+ \left[-\tfrac{r}{3}, \tfrac{r}{3}\right]^d$};

            \node[dot] at (2.5, 1.5) {};
            \node[cell] at (2.5, 1.5) {$\boldsymbol{y} + \tfrac{2r}{3}\boldsymbol{e}_1 + 0\boldsymbol{e}_2$ \\[3pt] $+ \left[-\tfrac{r}{3}, \tfrac{r}{3}\right]^d$};

            \node[dot] at (0.5, 0.5) {};
            \node[cell] at (0.5, 0.5) {$\boldsymbol{y} - \tfrac{2r}{3}\boldsymbol{e}_1 - \tfrac{2r}{3}\boldsymbol{e}_2$ \\[3pt] $+ \left[-\tfrac{r}{3}, \tfrac{r}{3}\right]^d$};

            \node[dot] at (1.5, 0.5) {};
            \node[cell] at (1.5, 0.5) {$\boldsymbol{y} + 0\boldsymbol{e}_1 - \tfrac{2r}{3}\boldsymbol{e}_2$ \\[3pt] $+ \left[-\tfrac{r}{3}, \tfrac{r}{3}\right]^d$};

            \node[dot] at (2.5, 0.5) {};
            \node[cell] at (2.5, 0.5) {$\boldsymbol{y} + \tfrac{2r}{3}\boldsymbol{e}_1 - \tfrac{2r}{3}\boldsymbol{e}_2$ \\[3pt] $+ \left[-\tfrac{r}{3}, \tfrac{r}{3}\right]^d$};

            \draw[dashed, thin, gray!70!black] (0.5, 0.5) -- (0.5, -0.22);
            \draw[dashed, thin, gray!70!black] (1.0, 0.0) -- (1.0, -0.22);
            \draw[<->, thin, gray!80!black] (0.5, -0.15) -- (1.0, -0.15) node[midway, below=1pt, black] {$\tfrac{r}{3}$};

            \draw[dashed, thin, gray!70!black] (1.5, 0.5) -- (1.5, -0.22);
            \draw[dashed, thin, gray!70!black] (3.0, 0.0) -- (3.0, -0.22);
            \draw[<->, thin, gray!80!black] (1.5, -0.15) -- (3.0, -0.15) node[midway, below=1pt, black] {$r$};

            \draw[dashed, thin, gray!70!black] (0.5, 1.5) -- (-0.22, 1.5);
            \draw[dashed, thin, gray!70!black] (0.5, 2.5) -- (-0.22, 2.5);
            \draw[<->, thin, gray!80!black] (-0.15, 1.5) -- (-0.15, 2.5) node[midway, left=2pt, black] {$\tfrac{2r}{3}$};
            }
        \end{tikzpicture}
        \caption{Division $D(Q)$ of a cube $Q = \boldsymbol{y} + [-r,r]^d$.}
        \label{fig:division}
    \end{figure}
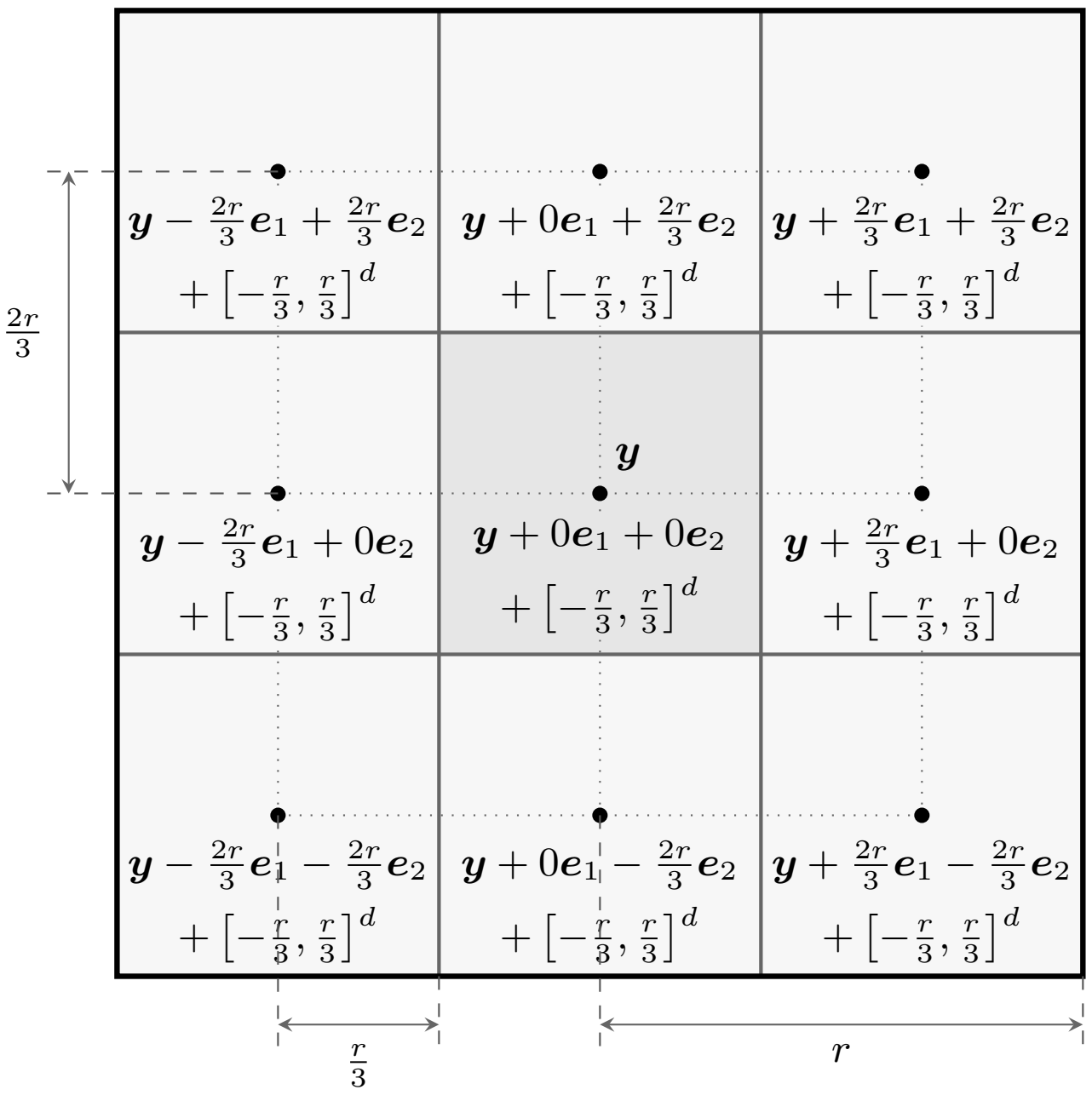
    
    We let
    \begin{equation*}
        D(\mathcal{Q}) = \bigcup_{Q\in \mathcal{Q}} D(Q)
    \end{equation*}
    and define a sequence $w_n = w_{D^n(\mathcal{Q})}$. Then summing up the above properties we have:
    \begin{itemize}
        \item $w_n(\boldsymbol{x}) > 0$ for all $n\in \mathbb{N}_+$ and a.e. $\boldsymbol{x}\in \Omega$,
        \item $\left|\nabla w_n(\boldsymbol{x})\right| = 1$ for all $n\in \mathbb{N}_+$ and a.e. $\boldsymbol{x}\in \Omega$,
        \item $\|w_n\|_\infty \leq \left(\frac{1}{3}\right)^n R_{\mathcal{Q}}$.
    \end{itemize}
    We show that $w_n\rightharpoonup 0$ in $W^{1,p}_0(\Omega)$. Indeed, the sequence $(w_n)$ is bounded in $W^{1,p}_0(\Omega)$ since ${\left\|w_n\right\|_{W^{1,p}_0} = |\Omega|^{\frac{1}{p}}}$. Moreover $w_n \to 0$ in $C\left(\overline{\Omega}\right)$. Consequently, each subsequence of $(w_n)$ has a~further subsequence weakly convergent to $0$. Hence the entire sequence $(w_n)$ converges weakly to $0$. Consequently, we can extract a subsequence $(u_n)$ of $(w_n)$ such that all assertions of the lemma are satisfied.
\end{proof}

\section{Properties of the energy functional}
\label{sec:energy_properties}

In this paper, we deal with equations on bounded domains, which is reflected in the following assumption:
\begin{equation}
    \label{ass:Omega}
    \tag{$\Omega$}
    \boxed{
    \text{$\Omega\subset \mathbb{R}^d$ is an open, bounded and non-empty set, $d \in \mathbb{N}_+$.}}
\end{equation}
Additionally, we adopt the assumption that is also made by other authors considering competing $(p,q)$-Laplacians:
\begin{equation}
    \label{ass:p<q}
    \tag{$p,q$}
    \boxed{
    \text{$1 < q < p < \infty$.}}
\end{equation}
Under \eqref{ass:Omega} and \eqref{ass:p<q}, we define the functional $\mathcal{D} \colon W^{1,p}_0(\Omega) \longrightarrow \mathbb{R}$ by the formula
\begin{equation}
    \label{eq:D_definition}
    \mathcal{D}(u) \coloneq \tfrac{1}{p} \int_{\Omega} |\nabla u(\boldsymbol{x})|^p \, d\boldsymbol{x} - \tfrac{1}{q} \int_{\Omega} |\nabla u(\boldsymbol{x})|^q \, d\boldsymbol{x} \quad \text{for every } u \in W^{1,p}_0(\Omega).
\end{equation}
In the preliminaries, we recalled typical assumptions from nonlinear functional analysis that lead to the existence of solutions for problems whose action functional is coercive, such as the functional $\mathcal{D}$ considered here. To better see its coercivity, it suffices to write it in the form
\begin{equation*}
    \mathcal{D}(u) = \int_\Omega \Phi\bigl(|\nabla u(\boldsymbol{x})|\bigr)\,d\boldsymbol{x}\qquad \text{for every } u \in W^{1,p}_0(\Omega),
\end{equation*}
where $\Phi \colon \mathbb{R}\longrightarrow \mathbb{R}$ is given by the formula
\begin{equation}
    \label{eq:Phi_definition}
    \Phi(t) = \tfrac{1}{p} |t|^p - \tfrac{1}{q} |t|^q\qquad \text{for every }t\in \mathbb{R}.
\end{equation}

\begin{figure}[htbp]
    \centering
    \begin{tikzpicture}[xscale=1.8, yscale=3.2, >=stealth]
        \draw[->, thin, gray!80!black] (-2.0, 0) -- (2.0, 0) node[right, black] {$t$};
        \draw[->, thin, gray!80!black] (0, -0.45) -- (0, 0.75) node[above, black] {$\Phi(t)$};

        \draw[dashed, gray] (-1, 0) -- (-1, -0.25) -- (1, -0.25) -- (1, 0);

        \node[above] at (-1, 0) {$-1$};
        \node[above] at (1, 0) {$1$};
        \node[below=2pt, fill=white, inner sep=1.5pt, font=\footnotesize] at (-0.2, -0.25) {$\tfrac{q-p}{pq}$};

        \draw[thick, smooth, samples=100, domain=-1.68:1.68] 
            plot (\x, {0.25*(\x)^4 - 0.5*(\x)^2});

        \node[circle, fill, inner sep=1.2pt] at (-1, -0.25) {};
        \node[circle, fill, inner sep=1.2pt] at (1, -0.25) {};
    \end{tikzpicture}
    \caption{A graph of the function $\Phi$ (plotted for $p=4$ and $q=2$).}
    \label{fig:Phi_graph}
\end{figure}
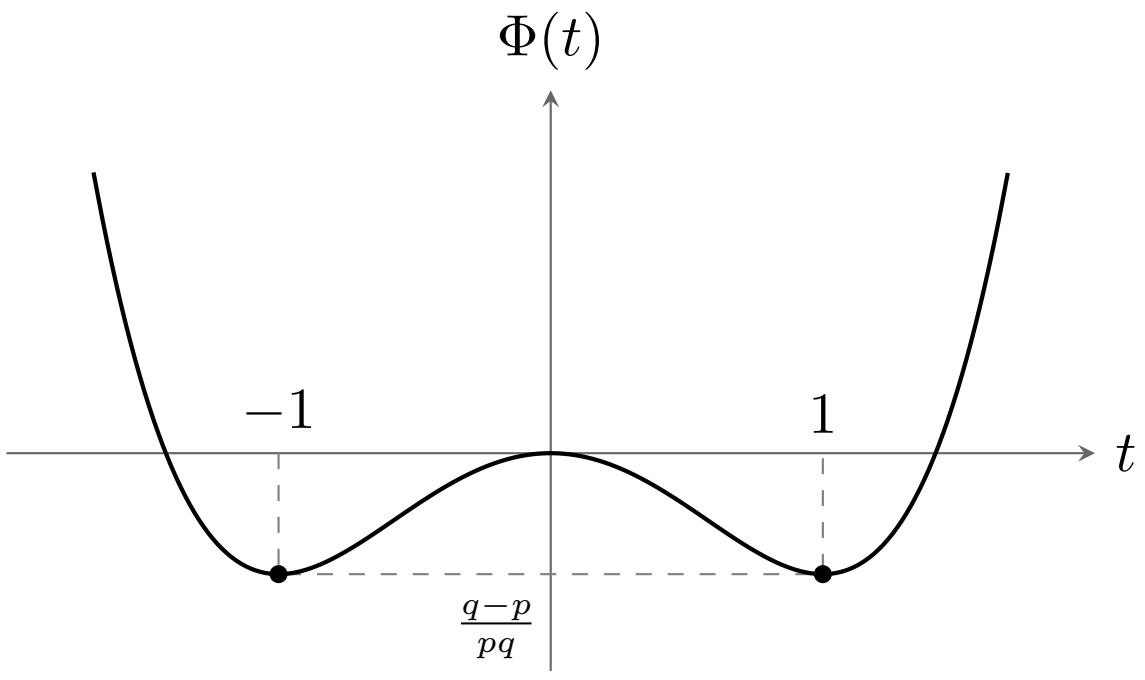

At the same time, the properties of the function $\Phi$ lead to a number of pathologies closely related to the fact that the global minimizers of $\mathcal{D}$ lie on the sphere of radius $|\Omega|^{1/p}$ in $W^{1,p}_0(\Omega)$. The above observations are described in detail in the following proposition.
\begin{proposition}
    \label{prop:D_properties}
    Assume that \eqref{ass:Omega} and \eqref{ass:p<q} hold. Define $\mathcal{D}$ by the formula \eqref{eq:D_definition}. Then $\mathcal{D}$ is neither w.l.s.c.\ nor s.w.l.s.c., and it does not satisfy (PS). Moreover, the operator $\mathcal{D}'$ is not pseudomonotone and does not satisfy condition (S)$_2$.
\end{proposition}
\begin{proof}
    Note that the function $\Phi$ given by the formula \eqref{eq:Phi_definition} achieves its global minimum at $t = 1$ and $t = -1$, with $\min_{t \in \mathbb{R}} \Phi(t) = \Phi(\pm 1) = \tfrac{q-p}{pq} < 0$. Hence
    \begin{equation*}
        \mathcal{D}(u)
        =
        \int_\Omega \Phi\bigl(|\nabla u(\boldsymbol{x})|\bigr)\,d\boldsymbol{x}
        \geq
        \tfrac{q - p}{pq}|\Omega|
        \qquad \text{for every } u \in W^{1,p}_0(\Omega).
    \end{equation*}
    Let $(u_n)$ be the sequence from Lemma \ref{lem:unit_gradient_sequence}. Since $|\nabla u_n(\boldsymbol{x})| = 1$ for a.e.\ $\boldsymbol{x} \in \Omega$, it follows that
    \begin{equation}
        \label{eq:u_n_are_critical}
        \mathcal{D}(u_n) = \min_{u \in W^{1,p}_0(\Omega)} \mathcal{D}(u) = \tfrac{q - p}{pq}|\Omega|\qquad \text{for every }n\in \mathbb{N}_+.
    \end{equation}
    Moreover, for every $v\in W^{1,p}_0(\Omega)$,
    \begin{equation*}
        \langle \mathcal{D}'(u_n), v\rangle
        =
        \int_\Omega \Bigl(|\nabla u_n(\boldsymbol{x})|^{p-2}-|\nabla u_n(\boldsymbol{x})|^{q-2}\Bigr)\nabla u_n(\boldsymbol{x})\cdot \nabla v(\boldsymbol{x})\,d\boldsymbol{x}
        =
        0,
    \end{equation*}
    and therefore $\mathcal{D}'(u_n)=0$ for every $n\in \mathbb{N}_+$.

    Note that $\mathcal{D}'(0) = 0$ and $\mathcal{D}(0) = 0$. Since $u_n \rightharpoonup 0$ in $W^{1,p}_0(\Omega)$, combining this with \eqref{eq:u_n_are_critical} gives
    \begin{equation*}
        \liminf_{n \to \infty} \mathcal{D}(u_n) = \tfrac{q - p}{pq}|\Omega| < 0 = \mathcal{D}(0),
    \end{equation*}
    so one can see that $\mathcal{D}$ is not s.w.l.s.c., and hence not w.l.s.c.\ either (see \cite[Section 1]{DeFigueiredo1989}). Furthermore, $\mathcal{D}$ does not satisfy (PS). Indeed, the sequence $(u_n)$ satisfies
    \begin{equation*}
        \sup_{n \in \mathbb{N}_+} |\mathcal{D}(u_n)| = \tfrac{p - q}{pq}|\Omega| < \infty
        \qquad \text{and} \qquad
        \mathcal{D}'(u_n) \to 0,
    \end{equation*}
    while it has no strongly convergent subsequence, since $u_n \rightharpoonup 0$ and
    \begin{equation*}
        \|u_n\|_{W^{1,p}_0}^p = \int_\Omega |\nabla u_n(\boldsymbol{x})|^p\,d\boldsymbol{x} = |\Omega| > 0
        \qquad \text{for every }n\in \mathbb{N}_+.
    \end{equation*}

    Note that $\mathcal{D}'$ is bounded on bounded sets. Hence, supposing $\mathcal{D}'$ is pseudomonotone, condition (H4) from \cite[Theorem 41.8]{Zeidler1985} can be used to show that $\mathcal{D}$ is also w.l.s.c. That clearly contradicts the already proven first part of this proposition. Finally, note that $u_n \rightharpoonup 0$ and $\mathcal{D}'(u_n) = 0 \to \mathcal{D}'(0)$. Hence, supposing that $\mathcal{D}'$ satisfies (S)$_2$, it follows that $u_n \to 0$ strongly in $W^{1,p}_0(\Omega)$, and consequently $\sqrt[p]{|\Omega|} = \|u_n\|_{W^{1,p}_0} \to 0$, whereby the desired contradiction is obtained.
\end{proof}

In order to handle the perturbation term in dimension $d$, we now impose an additional assumption on the nonlinear term $f \colon \Omega \times \mathbb{R} \longrightarrow \mathbb{R}$. Recall that $f$ is called \emph{a Carath{\'e}odory function} if $f(\boldsymbol{x}, \cdot)$ is continuous for a.e. $\boldsymbol{x} \in \Omega$ and $f(\boldsymbol{\cdot}, u)$ is measurable for every $u \in \mathbb{R}$. Depending on $d$ we impose the following assumption on $f$:
\begin{equation}
    \label{ass:f}
    \tag{$f$}
    \boxed{
    \begin{array}{l}
        \text{A Carath{\'e}odory function $f \colon \Omega \times \mathbb{R} \longrightarrow \mathbb{R}$ satisfies the growth condition:}\\[1ex]
        \quad\bullet\ \text{if $d < p$, for each $R > 0$ there exists $f_R \in L^1(\Omega)$ such that:}\\
        \qquad\quad |f(\boldsymbol{x}, u)| \le f_R(\boldsymbol{x}) \quad \text{for all } |u| \le R\text{ and a.e. $\boldsymbol{x} \in \Omega$};\\[1.5ex]
        \quad\bullet\ \text{if $d = p$, there exist $s > 1$, $f_0 \in L^s(\Omega)$, $C > 0$, and $1 \le r < \infty$ such that:}\\
        \qquad\quad |f(\boldsymbol{x}, u)| \le f_0(\boldsymbol{x}) + C |u|^{r-1} \quad \text{for all } u \in \mathbb{R}\text{ and a.e. $\boldsymbol{x} \in \Omega$};\\[1.5ex]
        \quad\bullet\ \text{if $d > p$, there exist $f_0 \in L^{(p^*)'}(\Omega)$, $C > 0$, and $1 \le r < p^*$ such that:}\\
        \qquad\quad |f(\boldsymbol{x}, u)| \le f_0(\boldsymbol{x}) + C |u|^{r-1} \quad \text{for all } u \in \mathbb{R}\text{ and a.e. $\boldsymbol{x} \in \Omega$}.
    \end{array}
    }
\end{equation}
Here $p^* = \frac{dp}{d-p}$ denotes the critical Sobolev exponent and $(p^*)' = \frac{dp}{dp-d+p}$ stands for its Hölder conjugate.

The antiderivative $F\colon \Omega \times \mathbb{R}\longrightarrow \mathbb{R}$ (with respect to the second variable) of the function $f$ is defined as
\begin{equation}
    \label{eq:F_definition}
    F(\boldsymbol{x}, u) \coloneq \int_0^u f(\boldsymbol{x}, v) \, dv \quad \text{for all } u \in \mathbb{R} \text{ and a.e.\ } \boldsymbol{x} \in \Omega.
\end{equation}

\begin{proposition}
    \label{prop:F_properties}
    Assume that \eqref{ass:f} holds. Then the functional $\mathcal{F} \colon W^{1,p}_0(\Omega) \longrightarrow \mathbb{R}$ defined by
    \begin{equation}
        \label{eq:F_functional}
        \mathcal{F}(u) \coloneq \int_{\Omega} F\big(\boldsymbol{x}, u(\boldsymbol{x})\big) \, d\boldsymbol{x} \qquad \text{for every } u \in W^{1,p}_0(\Omega),
    \end{equation}
    with $F$ given by \eqref{eq:F_definition}, is sequentially weakly continuous and of class $C^1$. Its derivative ${\mathcal{F}' \colon W^{1,p}_0(\Omega) \longrightarrow W^{-1,p'}(\Omega)}$ reads
    \begin{equation*}
        \langle \mathcal{F}'(u), v \rangle = \int_{\Omega} f\big(\boldsymbol{x}, u(\boldsymbol{x})\big) v(\boldsymbol{x}) \, d\boldsymbol{x} \qquad \text{for all } u, v \in W^{1,p}_0(\Omega)
    \end{equation*}
    and it is \emph{a strongly continuous operator}, that is, ${u_n \rightharpoonup u}$ in $W^{1,p}_0(\Omega)$ implies ${\mathcal{F}'(u_n) \to \mathcal{F}'(u)}$ in $W^{-1,p'}(\Omega)$.
\end{proposition}
\begin{proof}
    If $d > p$, we take any $a \in (1,p^*)$ satisfying $a \geq (r-1)\left(p^*\right)'$. Then by the assumption \eqref{ass:f} and \cite[Theorem~2.3 and Theorem~2.8]{DeFigueiredo1989}, the mapping
    \begin{equation*}
        L^{a}(\Omega) \ni u \longmapsto f\big(\boldsymbol{\cdot}, u(\boldsymbol{\cdot})\big) \in L^{(p^*)'}(\Omega)
    \end{equation*}
    is continuous. Since the Sobolev embedding $W^{1,p}_0(\Omega) \hookrightarrow L^{p^*}(\Omega)$ is continuous and the Rellich--Kondrachov theorem (see \cite[Theorem~7.22]{GilbargTrudinger2001}) ensures that the embeddings $W^{1,p}_0(\Omega) \hookrightarrow L^{a}(\Omega)$ and $W^{1,p}_0(\Omega) \hookrightarrow L^{r}(\Omega)$ are compact, the strong continuity of $\mathcal{F}'$, the sequential weak continuity of $\mathcal{F}$ and the regularity $\mathcal{F} \in C^1\big(W^{1,p}_0(\Omega)\big)$ follow.

    If $d = p$, take any finite $a > 1$ satisfying $a \geq (r-1)s$. Then, by \cite[Theorem~2.3 and Theorem~2.8]{DeFigueiredo1989} and by the assumption \eqref{ass:f}, the mapping
    \begin{equation*}
        L^{a}(\Omega) \ni u \longmapsto f\big(\boldsymbol{\cdot}, u(\boldsymbol{\cdot})\big) \in L^{s}(\Omega)
    \end{equation*}
    is continuous. Since the Rellich--Kondrachov theorem gives that the embeddings $W^{1,p}_0(\Omega) \hookrightarrow L^{a}(\Omega)$, $W^{1,p}_0(\Omega) \hookrightarrow L^{s'}(\Omega)$ and $W^{1,p}_0(\Omega) \hookrightarrow L^{r}(\Omega)$ are compact, all the desired properties of $\mathcal{F}$ follow as above.
    
    If $d < p$, then the embedding $W^{1,p}_0(\Omega) \hookrightarrow C(\overline{\Omega})$ is compact (see again \cite[Theorem~7.22]{GilbargTrudinger2001}). Using assumption \eqref{ass:f} and the Lebesgue dominated convergence theorem, one easily deduces the sequential weak continuity of $\mathcal{F}$, its G\^ateaux differentiability, and the strong continuity of $\mathcal{F}'$, which in turn yields $\mathcal{F} \in C^1\big(W^{1,p}_0(\Omega)\big)$ by \cite[Proposition~8.6]{DeFigueiredo1989}. The proof is standard, and hence we omit the details.
\end{proof}
Weak solutions to \eqref{eq:main_problem} coincide with critical points of the energy functional ${\mathcal{E} \colon W^{1,p}_0(\Omega) \longrightarrow \mathbb{R}}$ given by
\begin{equation}
    \label{eq:energy_functional}
    \mathcal{E}(u) \coloneq \tfrac{1}{p} \int_\Omega |\nabla u(\boldsymbol{x})|^p  \, d\boldsymbol{x} - \tfrac{1}{q} \int_\Omega |\nabla u(\boldsymbol{x})|^q \, d\boldsymbol{x} - \int_\Omega F\big(\boldsymbol{x}, u(\boldsymbol{x})\big) \, d\boldsymbol{x} \quad \text{for every } u \in W^{1,p}_0(\Omega),
\end{equation}
where $F$ is defined in \eqref{eq:F_definition}. Note that $\mathcal{E}$ given by \eqref{eq:energy_functional} can be expressed as $\mathcal{E} = \mathcal{D} - \mathcal{F}$. Consequently, the following is obtained.

\begin{proposition}
    \label{prop:E_properties}
    Let \eqref{ass:Omega}, \eqref{ass:p<q} and \eqref{ass:f} hold. Define the functional $\mathcal{E} \colon W^{1,p}_0(\Omega) \longrightarrow \mathbb{R}$ by formula \eqref{eq:energy_functional} and assume additionally that $f(\cdot, 0) \equiv 0$. Then $\mathcal{E}$ is neither w.l.s.c.\ nor s.w.l.s.c., and it does not satisfy (PS). Moreover, $\mathcal{E}'$ is neither pseudomonotone nor satisfies condition (S)$_2$.
\end{proposition}
\begin{proof}
    Supposing that $\mathcal{E}$ is s.w.l.s.c.\ (or w.l.s.c.), since $\mathcal{F}$ is sequentially weakly continuous by Proposition \ref{prop:F_properties}, the sum $\mathcal{D} = \mathcal{E} + \mathcal{F}$ would also be s.w.l.s.c., which contradicts Proposition \ref{prop:D_properties}. To~demonstrate that $\mathcal{E}$ does not satisfy (PS), the same sequence $(u_n)$ as in the proof of Proposition \ref{prop:D_properties} can be considered. Indeed, since ${u_n \rightharpoonup 0}$, we have ${\mathcal{F}(u_n) \to \mathcal{F}(0) = 0}$ and ${\mathcal{F}'(u_n) \to \mathcal{F}'(0) = 0}$ by Proposition \ref{prop:F_properties}. Hence, ${\mathcal{E}'(u_n) = \mathcal{D}'(u_n) - \mathcal{F}'(u_n) \to 0}$ and $(\mathcal{E}(u_n))$ is bounded, while $(u_n)$ has no convergent subsequence. Finally, since $\mathcal{F}'$ is strongly continuous, adding $\mathcal{F}'$ preserves pseudomonotonicity and condition (S)$_2$; thus, following the lines of the proof of Proposition \ref{prop:D_properties} for the sequence $(u_n)$ shows that $\mathcal{E}'$ is neither pseudomonotone nor satisfies condition (S)$_2$.
\end{proof}

The assumption $f(\cdot, 0) \equiv 0$ is a basic condition concerning the existence of nonzero critical points for semi-linear and quasilinear differential equations; see, for instance, \cite[Section 6.4]{Drabek2013}. It is a direct consequence of a standard assumption: ${f(\boldsymbol{x}, u) = o(|u|^{p-1})}$ (or ${f(\boldsymbol{x}, u) = o(|u|)}$) as $|u| \to 0$ uniformly for $\boldsymbol{x} \in \Omega$.

\section{Existence of a nonzero generalized solution}
\label{sec:generalized_solution}

The difficulties described in the previous section motivate a new definition, which is based on the Galerkin approximation method; see \cite{Motreanu2020} for a detailed description.

\begin{definition}[A generalized solution]
    A function $u \in W^{1,p}_0(\Omega)$ is called \emph{a generalized solution} to problem \eqref{eq:main_problem} if there exists a sequence $(u_n) \subset W^{1,p}_0(\Omega)$ such that $u_n \rightharpoonup u$ in $W^{1,p}_0(\Omega)$, $\langle \mathcal{E}'(u_n), v \rangle \to 0$ for each $v \in W^{1,p}_0(\Omega)$, and $\langle \mathcal{E}'(u_n), u_n - u \rangle \to 0$ as $n \to \infty$, where $\mathcal{E}$ is given by \eqref{eq:energy_functional}.
\end{definition}

It is important to note that multiple conceptions of generalized solutions currently exist. To the best of the author’s knowledge, the different notions of generalized solution used for competing operators have not yet been systematically compared. The differentiability of the functional $\mathcal{E}$ was a motivation to adapt the definition provided in \cite{Diblik2025} for the single-valued operator. So far, all results concerning the existence of a generalized solution have addressed only the existence itself. To the best of the author's knowledge, the existence of nonzero solutions in the generalized sense has not been analyzed. Note that in the case of the classical $p$-Laplacian, the existence of such solutions to problems of the type
\begin{equation}
    \label{eq:standard_pLaplacian}
    \begin{cases}
        -\Delta_p u(\boldsymbol{x}) = |u(\boldsymbol{x})|^{r-2}u(\boldsymbol{x}) & \text{for } \boldsymbol{x} \in \Omega, \\
        u(\boldsymbol{x}) = 0 & \text{for }\boldsymbol{x}\in \partial \Omega,
    \end{cases}
\end{equation}
for $r\in (1,p)$, can be deduced directly from the minimization of the associated energy functional. Indeed, taking $\mathcal{I}\colon W^{1,p}_0(\Omega) \longrightarrow \mathbb{R}$ given by the formula
\begin{equation*}
    \mathcal{I}(u) = \tfrac{1}{p}\int_\Omega |\nabla u(\boldsymbol{x})|^p \, d\boldsymbol{x} - \tfrac{1}{r}\int_{\Omega} |u(\boldsymbol{x})|^r \, d\boldsymbol{x}\quad \text{for every }u\in W^{1,p}_0(\Omega),
\end{equation*}
we immediately obtain
\begin{equation*}
    \inf_{u\in W^{1,p}_0(\Omega)}\mathcal{I}(u) < \mathcal{I}(0) = 0.
\end{equation*}
It follows from the fact that for a fixed $u\neq 0$, the function $t \longmapsto \mathcal{I}(tu)$ has a geometry similar to that of the function $\Phi$ (the details of such reasoning are well described in \cite[Sections~6.5 and~7.2]{Drabek2013}). For the functional $\mathcal{I}$, one can even show that it possesses infinitely many pairs of critical points,\linebreak see~\cite[Section~5.4]{Papageorgiou2019}. Standard arguments, however, require weak lower semicontinuity or the Palais--Smale condition. Both of these, as shown in the previous section, do not hold in the case of the functional $\mathcal{E}$. It seems all the more interesting that a case analogous to \eqref{eq:standard_pLaplacian} can be treated using the following theorem.
\begin{theorem}
    \label{thm:generalized_existence}
    Let \eqref{ass:Omega}, \eqref{ass:p<q} and \eqref{ass:f} hold and assume additionally that:
    \begin{itemize}
        \item $f(\cdot, 0) \equiv 0$,
        \item there exists $\varepsilon > 0$ such that $F(\boldsymbol{x}, u) > 0$ for all $u \in (0,\varepsilon)$ and a.e. $\boldsymbol{x}\in \Omega$,
        \item there exist  $\alpha \in\left(0, \tfrac{C_p}{p}\right)$ and $\beta \in L^1\big(\Omega, \mathbb{R}_+\big)$ such that 
        \begin{equation}
            \label{eq:F_growth}
            F(\boldsymbol{x}, u) \leq \alpha |u|^{p} + \beta(\boldsymbol{x}) \quad \text{for all }u\in \mathbb{R}\text{ and a.e. }\boldsymbol{x}\in \Omega.
        \end{equation}
    \end{itemize}
    Here $F$ is given by \eqref{eq:F_definition} and $C_p$ stands for the optimal constant in the Poincar{\'e} inequality \eqref{eq:Poincare}. Then problem \eqref{eq:main_problem} has a nonzero generalized solution $u$ (namely, the solution satisfies $u\not\equiv 0$).
\end{theorem}
\begin{proof}
    Note that the energy functional $\mathcal{E}$ given by \eqref{eq:energy_functional}, associated with problem \eqref{eq:main_problem}, is coercive. Indeed, since 
    \begin{equation*}
        \lim_{|t|\to \infty}\frac{\Phi(t)}{|t|^p} = \tfrac{1}{p}
    \end{equation*}
    for $\Phi$ given by \eqref{eq:Phi_definition}, there exists $c > 0$ such that
    \begin{equation*}
        \Phi(t) \geq \tfrac{C_p + p\alpha}{2p C_p}|t|^p - c\quad \text{for every }t\in \mathbb{R}.
    \end{equation*}
    Hence, using \eqref{eq:F_growth} and the Poincar{\'e} inequality \eqref{eq:Poincare}, we obtain
    \begin{equation}
        \label{eq:E_coercivity}
        \begin{split}
            \mathcal{E}(u)
            &=
            \int_\Omega \Phi\big( |\nabla u(\boldsymbol{x})|\big) \, d\boldsymbol{x}
            -
            \int_\Omega F\big(\boldsymbol{x}, u(\boldsymbol{x})\big) \, d\boldsymbol{x}\\
            &\geq
            \tfrac{C_p + p\alpha}{2p C_p} \|u\|_{W^{1,p}_0}^p
            -c|\Omega|
            -
            \alpha \int_\Omega |u(\boldsymbol{x})|^p \, d\boldsymbol{x}
            -
            \int_\Omega \beta(\boldsymbol{x}) \, d\boldsymbol{x}\\
            &\geq
            \tfrac{C_p - p\alpha}{2p C_p}\|u\|_{W^{1,p}_0}^p
            -
            c|\Omega|
            -
            \|\beta\|_{L^1(\Omega)}.
        \end{split}
    \end{equation}
    Therefore $\mathcal{E}$ is coercive and bounded from below. Consequently, by the Ekeland Variational Principle, there exists a sequence $(w_n) \subset W^{1,p}_0(\Omega)$ satisfying
    \begin{equation*}
        \mathcal{E}(w_n) \to \inf_{u \in W^{1,p}_0(\Omega)} \mathcal{E}(u) \quad \text{and} \quad \mathcal{E}'(w_n) \to 0 \quad \text{in } W^{-1,p'}(\Omega).
    \end{equation*}
    Such a sequence must be bounded in $W^{1,p}_0(\Omega)$ due to the coercivity of $\mathcal{E}$. Denote by $(u_n)$ a subsequence of $(w_n)$ that satisfies $u_n \rightharpoonup u_0$ in $W^{1,p}_0(\Omega)$ for some $u_0 \in W^{1,p}_0(\Omega)$ (such a subsequence exists by the reflexivity of $W^{1,p}_0(\Omega)$). Since $\|\mathcal{E}'(u_n)\|_{W^{-1,p'}} \to 0$ and $\sup_{n \in \mathbb{N}_+} \|u_n - u_0\|_{W^{1,p}_0} < \infty$, we have $\langle \mathcal{E}'(u_n), v \rangle \to 0$ for every $v \in W^{1,p}_0(\Omega)$ as well as
    \begin{equation*}
        \big|\langle \mathcal{E}'(u_n), u_n - u_0 \rangle\big| \leq \|\mathcal{E}'(u_n)\|_{W^{-1,p'}} \|u_n - u_0\|_{W^{1,p}_0} \to 0 \qquad \text{as } n \to \infty.
    \end{equation*}
    Thus $u_0$ is a generalized solution to \eqref{eq:main_problem}.

    Suppose now that $u_0 = 0$. Recall that $\mathcal{E} = \mathcal{D} - \mathcal{F}$, where $\mathcal{D}$ and $\mathcal{F}$ are defined by \eqref{eq:D_definition} and \eqref{eq:F_functional}, respectively. Since $u_n \rightharpoonup 0$, Proposition \ref{prop:F_properties} yields
    \begin{equation*}
        \mathcal{F}(u_n) \to \mathcal{F}(0) = 0.
    \end{equation*}
    On the other hand, Lemma \ref{lem:unit_gradient_sequence} allows us to take $u\in W^{1,p}_0(\Omega)$ satisfying:
    \begin{equation*}
        0 < u(\boldsymbol{x}) < \tfrac{\varepsilon}{2}\qquad \text{and} \qquad |\nabla u(\boldsymbol{x})| = 1\qquad \text{ for a.e. }\boldsymbol{x}\in \Omega.
    \end{equation*}
    Then
    \begin{equation}
        \label{eq:upper_inf}
        \inf_{v\in W^{1,p}_0(\Omega)} \mathcal{E}(v) \leq \mathcal{E}(u) = \mathcal{D}(u) - \mathcal{F}(u) < \tfrac{q-p}{pq}|\Omega|.
    \end{equation}
    However, the lower bound established in the proof of Proposition \ref{prop:D_properties} gives
    \begin{equation*}
        \mathcal{D}(u_n)\geq \tfrac{q-p}{pq}|\Omega| \qquad \text{for every }n\in \mathbb{N}_+,
    \end{equation*}
    and therefore, since $(u_n)$ is a minimizing sequence for $\mathcal{E}$, we get
    \begin{equation*}
        \inf_{u\in W^{1,p}_0(\Omega)} \mathcal{E}(u)
        =
        \lim_{n\to\infty}\mathcal{E}(u_n)
        =
        \lim_{n\to\infty}\bigl(\mathcal{D}(u_n)-\mathcal{F}(u_n)\bigr)
        \geq
        \tfrac{q-p}{pq}|\Omega|.
    \end{equation*}
    This contradicts \eqref{eq:upper_inf} and hence shows that $u_0 \neq 0$.
\end{proof}

Observe that the function $f\colon \Omega \times \mathbb{R}\longrightarrow \mathbb{R}$ defined by the formula
\begin{equation*}
    f(\boldsymbol{x}, u) = |u|^{r-2}u \quad \text{for all }u\in \mathbb{R} \text{ and every }\boldsymbol{x}\in \Omega,
\end{equation*}
with $r\in (1,p)$, satisfies all assumptions of Theorem \ref{thm:generalized_existence}. Consequently we can conclude that problem \eqref{eq:standard_pLaplacian} with $-\Delta_p$ replaced by the competing $(p,q)$-Laplacian has a nonzero generalized solution.
\begin{remark}
    \label{remark1}
    Let $\Omega\subset \mathbb{R}^d$, $d\in \mathbb{N}_+$, be open and bounded. Take $p > 1$ and $q,r\in (1,p)$. Then the problem 
    \begin{equation*}
        \begin{cases}
            -\Delta_p u(\boldsymbol{x}) + \Delta_q u(\boldsymbol{x}) = |u(\boldsymbol{x})|^{r-2}u(\boldsymbol{x}) & \text{for } \boldsymbol{x} \in \Omega, \\
            u(\boldsymbol{x}) = 0 & \text{for }\boldsymbol{x}\in \partial \Omega,
        \end{cases}
    \end{equation*}
    has at least one nonzero generalized solution.
\end{remark}

\section{Existence of a nonzero weak solution on an interval}
\label{sec:weak_solution}

Considering problem \eqref{eq:main_problem} in $\Omega = (-1,1)$ leads to its one-dimensional formulation \eqref{eq:main_problem-1dim}. The energy $\mathcal{E}\colon W^{1,p}_0(-1,1)\longrightarrow \mathbb{R}$ defined by \eqref{eq:energy_functional}, in this particular case, reads
\begin{equation}
    \label{eq:energy_functional-1dim}
    \mathcal{E}(u) = \int_{-1}^1 \Phi\big(u'(x)\big)\, dx - \int_{-1}^1 F\big(x, u(x)\big)\, dx.
\end{equation}
The functions $\Phi$ and $F$ are defined respectively by \eqref{eq:Phi_definition} and \eqref{eq:F_definition}. To properly define $\mathcal{E}$ we assume that the nonlinear term $f\colon (-1,1)\times \mathbb{R}\longrightarrow \mathbb{R}$ satisfies assumption \eqref{ass:f}, which has the following form:
\begin{equation}
    \label{ass:f1}
    \tag{$f_1$}
    \boxed{
    \begin{array}{l}
        \text{The function $f \colon (-1,1) \times \mathbb{R} \longrightarrow \mathbb{R}$ is a Carath{\'e}odory function.}\\
        \text{Moreover, for each $R > 0$ there exists $f_R \in L^1(-1,1)$ such that:}\\
        \qquad\quad |f(x, u)| \le f_R(x) \quad \text{for all } |u| \le R\text{ and a.e. $x \in (-1,1)$}.
    \end{array}
    }
\end{equation}
We impose additional assumptions on the sign and symmetry of the nonlinearity $f$ of the form
\begin{equation}
    \label{ass:fsym}
    \tag{$f_{\text{sym}}$}
    \boxed{
    f(x,u) = f(-x,u) = -f(x,-u)\quad \text{for all }u\in \mathbb{R}\text{ and a.e. }x\in (-1, 1).
    }
\end{equation}
Such an assumption allows us to study only symmetric solutions to \eqref{eq:main_problem-1dim}, namely those satisfying ${u(x) = u(-x)}$ for $x\in (-1,1)$, due to the Palais Principle, see \cite{Palais}. The principle itself, however, does not settle the question of the minimization of the functional. Therefore we define the functional $\mathcal{J} \colon \widetilde{W^{1,p}_0}(-1,0) \longrightarrow \mathbb{R}$ by the formula
\begin{equation}
    \label{eq:definition_J}
    \mathcal{J}(u) = \int_{-1}^0 \Phi\big(u'(x)\big)\, dx - \int_{-1}^0 F\big(x, u(x)\big)\, dx\quad \text{for every }u\in \widetilde{W^{1,p}_0}(-1,0),
\end{equation}
where 
\begin{equation*}
    \widetilde{W^{1,p}_0}(-1,0) = \left\{u\in W^{1,p}(-1,0) : u(-1) = 0\right\}
\end{equation*}
is equipped with the $W^{1,p}_0$-norm, that is
\begin{equation*}
    \|u\|_{W^{1,p}_0} = \left(\int_{-1}^{0} \left|u'(x)\right|^p dx \right)^\frac{1}{p}\quad \text{for every }u\in \widetilde{W^{1,p}_0}(-1,0).
\end{equation*}
The connection between the global minimizers of $\mathcal{E}$ and $\mathcal{J}$ is described in the following lemma.
\begin{lemma}
    \label{lem:comparison_of_minimizers}
    Assume that \eqref{ass:p<q}, \eqref{ass:f1} and \eqref{ass:fsym} hold. Define $\mathcal{E}$ and $\mathcal{J}$ by \eqref{eq:energy_functional-1dim} and \eqref{eq:definition_J}, respectively. Then:
    \begin{itemize}
        \item if $u\in W^{1,p}_0(-1,1)$ satisfies 
        \begin{equation}
            \label{eq:minimization_of_E}
            \mathcal{E}(u) = \min_{v\in W^{1,p}_0(-1,1)} \mathcal{E}(v),
        \end{equation}
        then, taking $u_1 = u|_{[-1,0]}$, we get
        \begin{equation}
            \label{eq:minimization_of_J}
            \mathcal{J}(u_1) = \min_{v\in \widetilde{W^{1,p}_0}(-1,0)} \mathcal{J}(v),
        \end{equation}
        \item on the other hand, if $u_1\in \widetilde{W^{1,p}_0}(-1,0)$ satisfies \eqref{eq:minimization_of_J}, then the symmetric extension $\widetilde{u_1}$ of $u_1$ defined as
        \begin{equation}
            \label{eq:symmetric_extension}
            \widetilde{u_1}(x) = u_1(-|x|)\quad \text{for every }x\in [-1,1]
        \end{equation}
        satisfies \eqref{eq:minimization_of_E}.
    \end{itemize}
\end{lemma}
\begin{proof}
    Let us first note that as a consequence of \eqref{ass:fsym} we get
    \begin{equation*}
        F(x,u) = F(x, -u) = F(-x,u)\quad \text{for all }u\in \mathbb{R}\text{ and a.e. }x\in (-1, 1).
    \end{equation*}
    Take any $v\in W^{1,p}_0(-1,1)$. Set $v_1 = v|_{[-1,0]}$ and define $v_2\in \widetilde{W^{1,p}_0}(-1,0)$ by the formula
    \begin{equation}
        \label{eq:v2}
        v_2(x) = v(-x) \quad \text{for every }x\in [-1,0].
    \end{equation}
    Then, $v_2'(x) = -v'(-x)$ for a.e. $x\in (-1,0)$. Moreover, by the evenness of $\Phi$ and the symmetry of $F$ we have
    \begin{equation}
        \label{eq:E_J_decomposition}
        \begin{split}
            \mathcal{E}(v) & = \int_{-1}^0 \Phi\big(v'(x)\big)\,dx + \int_0^1 \Phi\big(v'(x)\big)\,dx - \int_{-1}^0 F\big(x, v(x)\big)\,dx - \int_0^1 F\big(x, v(x)\big)\,dx \\
            & = \int_{-1}^0 \Phi\big(v'(x)\big)\,dx - \int_{-1}^0 F\big(x, v(x)\big)\,dx + \int_{-1}^0 \Phi\big(-v'(-x)\big)\,dx - \int_{-1}^0 F\big(x, v(-x)\big)\,dx \\
            & = \mathcal{J}(v_1) + \mathcal{J}(v_2).
        \end{split}
    \end{equation}
    In particular, if ${\widetilde{w}\in W^{1,p}_0(-1,1)}$ is the symmetric extension of some ${w\in \widetilde{W^{1,p}_0}(-1,0)}$, then ${v_1 = v_2 = w}$, and hence \eqref{eq:E_J_decomposition} yields
    \begin{equation}
        \label{eq:E_J_symmetric}
        \mathcal{E}(\widetilde{w}) = 2\mathcal{J}(w).
    \end{equation}

    Assume first that \eqref{eq:minimization_of_E} holds for some $u\in W^{1,p}_0(-1,1)$. Set $u_1 = u|_{[-1,0]}$ and $u_2(x) = u(-x)$ for $x\in [-1,0]$ and consider the symmetric extension $\widetilde{u_2}$ of $u_2$ defined as in \eqref{eq:symmetric_extension}. From \eqref{eq:minimization_of_E}, \eqref{eq:E_J_decomposition} and \eqref{eq:E_J_symmetric} we immediately deduce that
    \begin{equation*}
        \mathcal{J}(u_1) + \mathcal{J}(u_2) = \mathcal{E}(u) \leq \mathcal{E}(\widetilde{u_2}) = 2\mathcal{J}(u_2),
    \end{equation*}
    which implies $\mathcal{J}(u_1) \leq \mathcal{J}(u_2)$, and consequently $2\mathcal{J}(u_1) \leq \mathcal{E}(u)$. Hence, taking any ${w \in \widetilde{W^{1,p}_0}(-1,0)}$ and defining its symmetric extension $\widetilde{w}$ by \eqref{eq:symmetric_extension} we can use \eqref{eq:minimization_of_E} together with \eqref{eq:E_J_symmetric} to get
    \begin{equation*}
        \mathcal{J}(u_1) \leq \tfrac{1}{2}\mathcal{E}(u) \leq \tfrac{1}{2}\mathcal{E}(\widetilde{w}) = \mathcal{J}(w).
    \end{equation*}
    Therefore \eqref{eq:minimization_of_J} holds.

    Now, assume that \eqref{eq:minimization_of_J} is satisfied for some $u_1\in \widetilde{W^{1,p}_0}(-1,0)$ and define $u\in W^{1,p}_0(-1,1)$ as the symmetric extension of $u_1$, that is ${u = \widetilde{u_1}}$. Take any ${v\in W^{1,p}_0(-1,1)}$ and define ${v_1, v_2 \in \widetilde{W^{1,p}_0}(-1,0)}$, respectively, by ${v_1 = v|_{[-1,0]}}$ and by \eqref{eq:v2}. Then, by \eqref{eq:minimization_of_J}, \eqref{eq:E_J_decomposition} and \eqref{eq:E_J_symmetric}, we have
    \begin{equation*}
        \mathcal{E}(v) = \mathcal{J}(v_1) + \mathcal{J}(v_2) \geq 2\mathcal{J}(u_1) = \mathcal{E}(u),
    \end{equation*}
    which means that $u$ satisfies \eqref{eq:minimization_of_E}.
\end{proof}

To establish the existence of minimizers despite the pathological properties of the energy functional~$\mathcal{E}$ observed in Section \ref{sec:energy_properties}, we define, under \eqref{ass:p<q}, the convex envelope $\widehat{\Phi}\colon \mathbb{R}\longrightarrow \mathbb{R}$ of $\Phi$ as follows:
\begin{equation*}
    \widehat{\Phi}(t) = 
    \begin{cases}
        \frac{1}{p}|t|^p - \frac{1}{q}|t|^q & \text{if } |t| \geq 1,\\
        \frac{q - p}{pq} & \text{if } t \in (-1,1).
    \end{cases}
\end{equation*}

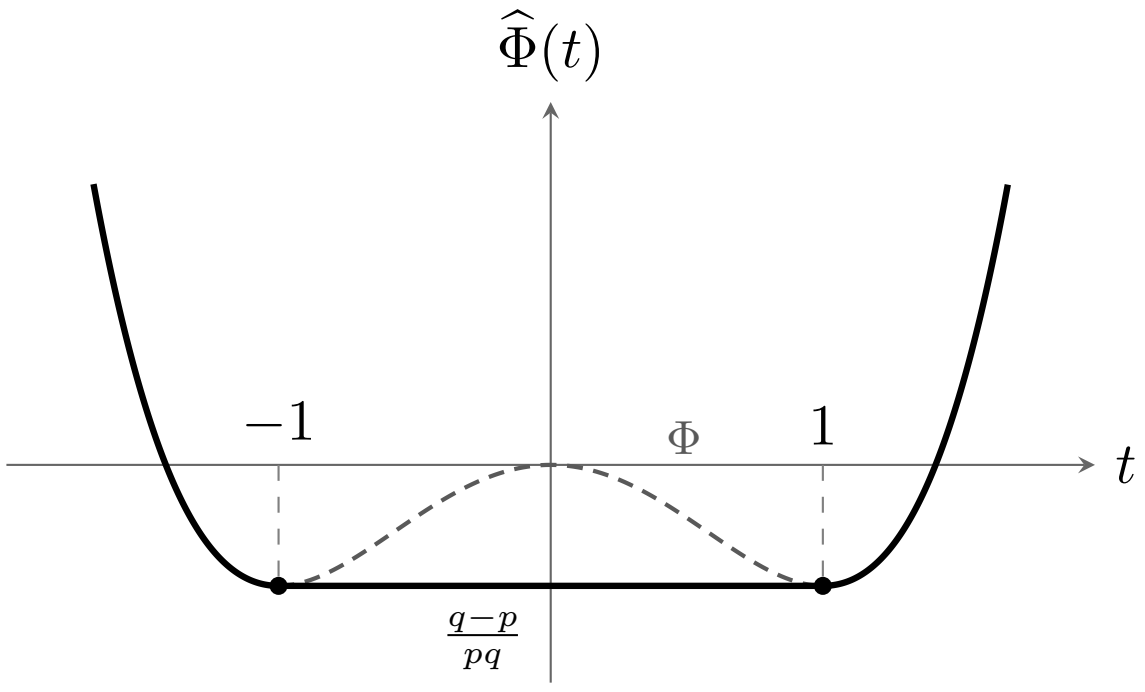
\begin{figure}[!htbp]
    \centering
    \begin{tikzpicture}[xscale=1.8, yscale=3.2, >=stealth]
        \draw[->, thin, gray!80!black] (-2.0, 0) -- (2.0, 0) node[right, black] {$t$};
        \draw[->, thin, gray!80!black] (0, -0.45) -- (0, 0.75) node[above, black] {$\widehat{\Phi}(t)$};

        \draw[dashed, gray] (-1, 0) -- (-1, -0.25);
        \draw[dashed, gray] (1, 0) -- (1, -0.25);

        \node[above] at (-1, 0) {$-1$};
        \node[above] at (1, 0) {$1$};
        \node[below=2pt, fill=white, inner sep=1.5pt, font=\footnotesize] at (-0.25, -0.25) {$\tfrac{q-p}{pq}$};

        \draw[thick, densely dashed, gray!70!black, smooth, samples=60, domain=-1:1] 
            plot (\x, {0.25*(\x)^4 - 0.5*(\x)^2});
        \node[above right, gray!70!black, font=\footnotesize] at (0.35, -0.02) {$\Phi$};

        \draw[very thick, smooth, samples=50, domain=-1.68:-1] 
            plot (\x, {0.25*(\x)^4 - 0.5*(\x)^2});
        \draw[very thick] (-1, -0.25) -- (1, -0.25) node[midway, below=2pt, right=8pt, font=\footnotesize]{};
        \draw[very thick, smooth, samples=50, domain=1:1.68] 
            plot (\x, {0.25*(\x)^4 - 0.5*(\x)^2});

        \node[circle, fill, inner sep=1.2pt] at (-1, -0.25) {};
        \node[circle, fill, inner sep=1.2pt] at (1, -0.25) {};
    \end{tikzpicture}
    \caption{A graph of the functions $\Phi$ and $\widehat{\Phi}$ (plotted for $p=4$ and $q=2$).}
    \label{fig:Phi_hat_graph}
\end{figure}

We now define $\widehat{\mathcal{J}}\colon \widetilde{W^{1,p}_0}(-1,0)\longrightarrow \mathbb{R}$ by the formula
\begin{equation}
    \label{eq:definition_Jhat}
    \widehat{\mathcal{J}}(u) = \int_{-1}^0 \widehat{\Phi}\big(u'(x)\big)\, dx - \int_{-1}^0 F\big(x, u(x)\big)\, dx\quad \text{for every }u\in \widetilde{W^{1,p}_0}(-1,0).
\end{equation}
It should be noted that the functional $\widehat{\mathcal{J}}$ does not necessarily coincide with the convex envelope of $\mathcal{J}$. While this choice of notation might appear somewhat misleading, we use it to emphasize the fact that the definition of $\widehat{\mathcal{J}}$ relies directly on $\widehat{\Phi}$. For a better control of minimizers of $\widehat{\mathcal{J}}$ we add assumptions
\begin{equation}
    \label{ass:f>=}
    \tag{$f_\geq$}
    \boxed{
    f(x,u) \geq 0\quad \text{for all }u\geq 0 \text{ and a.e. }x\in (-1, 1)
    }
\end{equation}
and 
\begin{equation}
    \label{ass:f>}
    \tag{$f_>$}
    \boxed{
    f(x,u) > 0\quad \text{for all }u > 0 \text{ and a.e. }x\in (-1, 1).
    }
\end{equation}
The consideration of envelopes, also known as functional relaxation, is a standard approach widely used in the literature, see for instance \cite{Celada1D, CeladaCoV}. The reasoning presented below is analogous to that presented in \cite{Kromer}, with the difference that the work of Kr\"omer investigates the multidimensional case for a~considerably more general functional $\Phi$, whereas the nonlinearity $f$ is independent of $x$. In~the final part of the paper, we will return to a comparison of the result below with \cite{Kromer}, based on the special case considered in Remark~\ref{remark1}.

Before moving to the main existence result, let us note that, by \cite{DrabekKrejciTakac1999}, in the one-dimensional case, the Poincar{\'e} constant~\eqref{eq:Poincare} in $W^{1,p}_0(-1,1)$ reads
\begin{equation}
    \label{eq:Poincare1D}
    C_p = (p - 1)\left(\frac{\pi}{p\sin\left(\frac{\pi}{p}\right)}\right)^p.
\end{equation}
\begin{theorem}
    Assume that \eqref{ass:p<q}, \eqref{ass:f1}, \eqref{ass:fsym} and \eqref{ass:f>=} hold. Then for every weak solution ${u \in W^{1,p}_0(-1,1)}$ to \eqref{eq:main_problem-1dim} being a global minimizer of $\mathcal{E}$, given by \eqref{eq:energy_functional-1dim}, there exists a weak solution $u^*\in W^{1,p}_0(-1,1)$ satisfying
    \begin{equation}
        \label{eq:properties_of_weak_solution}
        u^*\notin C^1(-1,1),\qquad  u^*(x) > 0 \text{ for a.e. }x\in (-1,1)\qquad\text{and}\qquad \mathcal{E}(u) = \mathcal{E}(u^*).
    \end{equation}
    The solution minimizing $\mathcal{E}$ exists provided there exist $\alpha < \frac{p - 1}{p}\left(\frac{\pi}{p\sin\left(\frac{\pi}{p}\right)}\right)^p$ and $\beta \in L^1\big((-1,1), \mathbb{R}_+\big)$ such that 
    \begin{equation}
        \label{eq:F_growth1D}
        F(x, u) \leq \alpha |u|^{p} + \beta(x) \quad \text{for all }u\in \mathbb{R}\text{ and a.e. }x\in (-1,1).
    \end{equation}
    Here $F$ is the antiderivative of $f$ given by \eqref{eq:F_definition}.

    If additionally \eqref{ass:f>} holds, then every weak solution $u$ to \eqref{eq:main_problem-1dim} that minimizes $\mathcal{E}$ is not of class $C^1$ and does not vanish in $(-1,1)$.
\end{theorem}
\begin{proof}
    Assume first that \eqref{ass:p<q}, \eqref{ass:f1}, \eqref{ass:fsym} and \eqref{ass:f>=} are satisfied. For any $v \in \widetilde{W^{1,p}_0}(-1,0)$, define $\overline{v} \in \widetilde{W^{1,p}_0}(-1,0)$ by the formula
    \begin{equation*}
        \overline{v}(x) = \int_{-1}^x \big(|v'(\xi)| \vee 1\big)\, d\xi \quad\text{for every }x\in [-1,0],
    \end{equation*}
    where, as usual, $a\vee b$ denotes the maximum of $a$ and $b$. Note that we immediately obtain
    \begin{equation}
        \label{eq:u'>=u'}
        \overline{v}'(x) = |v'(x)| \vee 1 \geq |v'(x)|\quad \text{for a.e. }x\in (-1,0).
    \end{equation}
    Hence, since $\widehat{\Phi}$ is even, coincides with $\Phi$ on $\mathbb{R}\setminus (-1,1)$, and satisfies $\widehat{\Phi}(t) = \Phi(1) \leq \Phi(t)$ for $t \in (-1,1)$, we get
    \begin{equation}
        \label{eq:iqualitywrtPhi}
        \Phi\big(\overline{v}'(x)\big) = \widehat{\Phi}\big(\overline{v}'(x)\big) = \widehat{\Phi}\big(v'(x)\big) \leq \Phi\big(v'(x)\big)\quad \text{for a.e. }x\in (-1,0).
    \end{equation}
    Moreover, by \eqref{eq:u'>=u'}, for each $x \in (-1,0)$ we have
    \begin{equation*}
        \overline{v}(x) - |v(x)| \geq \int_{-1}^x \big(\overline{v}'(\xi) - |v'(\xi)| \big)\, d\xi \geq 0.
    \end{equation*}
    The non-negativity of the integrand above gives the following dichotomy:
    \begin{equation}
        \label{eq:dichotomy}
        \text{Either $|v| = \overline{v}$ or there exists $x^*\in (-1,0)$ such that $\overline{v}(x) > |v(x)|$ for every $x \in (x^*,0)$.}
    \end{equation}
    Now, let us observe that assumption \eqref{ass:fsym} provides evenness of $F(x,\cdot)$ for a.e. $x\in (-1,1)$, while \eqref{ass:f>=} gives that $F(x,\cdot)$ is non-decreasing on $[0,\infty)$ for a.e. $x\in (-1,1)$. Combining this with \eqref{eq:u'>=u'}, \eqref{eq:iqualitywrtPhi} and \eqref{eq:dichotomy} shows that for every $v \in \widetilde{W^{1,p}_0}(-1,0)$,
    \begin{equation}
        \label{eq:J_Jhat_chain}
        \begin{split}
            \mathcal{J}(\overline{v}) = \widehat{\mathcal{J}}(\overline{v}) & = \int_{-1}^0 \widehat{\Phi}\big(\overline{v}'(x)\big)\, dx - \int_{-1}^0 F\big(x, \overline{v}(x)\big)\, dx \\
            & \leq \int_{-1}^0 \widehat{\Phi}\big(v'(x)\big)\, dx - \int_{-1}^0 F\big(x, v(x)\big)\, dx = \widehat{\mathcal{J}}(v) \leq \mathcal{J}(v).
        \end{split}
    \end{equation}

    Let $u \in W^{1,p}_0(-1,1)$ be a weak solution to \eqref{eq:main_problem-1dim} that globally minimizes $\mathcal{E}$. By Lemma~\ref{lem:comparison_of_minimizers}, its restriction $v = u|_{[-1,0]} \in \widetilde{W^{1,p}_0}(-1,0)$ satisfies
    \begin{equation}
        \label{eq:u1_minimizes_J}
        \mathcal{J}(v) = \min_{\varphi\in \widetilde{W^{1,p}_0}(-1,0)} \mathcal{J}(\varphi).
    \end{equation}
    In view of \eqref{eq:J_Jhat_chain} and \eqref{eq:u1_minimizes_J}, it follows that 
    \begin{equation*}
        \mathcal{J}(\overline{v}) = \mathcal{J}(v) = \min_{\varphi\in \widetilde{W^{1,p}_0}(-1,0)} \mathcal{J}(\varphi).
    \end{equation*}
    Let $u^* \in W^{1,p}_0(-1,1)$ be the symmetric extension of $\overline{v}$ defined as in \eqref{eq:symmetric_extension}. By Lemma~\ref{lem:comparison_of_minimizers}, $u^*$ globally minimizes $\mathcal{E}$ on $W^{1,p}_0(-1,1)$, so $\mathcal{E}(u) = \mathcal{E}(u^*)$ and $u^*$ is a weak solution to \eqref{eq:main_problem-1dim}. Furthermore, by \eqref{eq:u'>=u'}, we have 
    \begin{equation*}
        u^*(x) = \overline{v}(-|x|) \geq 1 - |x| > 0\quad \text{for every }x\in (-1,1).
    \end{equation*}
    Finally, since $u^*$ is the symmetric extension of $\overline{v}$, we have 
    \begin{equation*}
        (u^*)'(x) = \overline{v}'(x) \geq 1\text{ for a.e. }x\in (-1,0)\quad\text{and}\quad (u^*)'(x) = -\overline{v}'(-x) \leq -1 \text{ for a.e. }x \in (0,1).
    \end{equation*}
    Consequently, $(u^*)'$ cannot be continuous at $0$, which means that $u^* \notin C^1(-1,1)$ and proves \eqref{eq:properties_of_weak_solution}.

    Now, let us prove the existence of a global minimizer of $\mathcal{E}$ under assumption \eqref{eq:F_growth1D}. Consider the relaxed functional $\widehat{\mathcal{J}}$ defined by \eqref{eq:definition_Jhat}. The principal part of $\widehat{\mathcal{J}}$, i.e., 
    \begin{equation*}
        \widetilde{W^{1,p}_0}(-1,0) \ni v \longmapsto \int_{-1}^0 \widehat{\Phi}\big(v'(x)\big)\, dx \in \mathbb{R},
    \end{equation*}
    is convex and continuous and hence w.l.s.c., see \cite[Theorem~1.4]{DeFigueiredo1989}. Therefore, it is also s.w.l.s.c.,\linebreak see \cite[Section~1]{DeFigueiredo1989}. The second term of $\widehat{\mathcal{J}}$, namely
    \begin{equation*}
        \widetilde{W^{1,p}_0}(-1,0) \ni v \longmapsto \int_{-1}^0 F\big(x, v(x)\big)\, dx \in \mathbb{R},
    \end{equation*}
    is sequentially weakly continuous by the same argument as in the case ${d < p}$ of Proposition~\ref{prop:F_properties} (since the embedding ${\widetilde{W^{1,p}_0}(-1,0) \hookrightarrow C([-1,0])}$ is compact). Hence, $\widehat{\mathcal{J}}$ is s.w.l.s.c. as a sum of a s.w.l.s.c. functional and a weakly continuous one. Coercivity of $\widehat{\mathcal{J}}$ follows by direct calculations analogous to \eqref{eq:E_coercivity}, using \eqref{eq:F_growth1D} and the Poincar\'e inequality \eqref{eq:Poincare1D}. Thus, by the direct method of the calculus of variations,\linebreak see~\cite[Theorem~1.2]{DeFigueiredo1989}, there exists $w \in \widetilde{W^{1,p}_0}(-1,0)$ such that
    \begin{equation*}
        \widehat{\mathcal{J}}(w) = \min_{v\in \widetilde{W^{1,p}_0}(-1,0)} \widehat{\mathcal{J}}(v).
    \end{equation*}
    Applying \eqref{eq:J_Jhat_chain} to $w$, we obtain
    \begin{equation*}
        \mathcal{J}(\overline{w}) = \widehat{\mathcal{J}}(\overline{w}) \leq \widehat{\mathcal{J}}(w) \leq \widehat{\mathcal{J}}(v) \leq \mathcal{J}(v) \quad \text{for every } v \in \widetilde{W^{1,p}_0}(-1,0).
    \end{equation*}
    Hence $\overline{w}$ minimizes $\mathcal{J}$ on $\widetilde{W^{1,p}_0}(-1,0)$, and by Lemma~\ref{lem:comparison_of_minimizers} its symmetric extension $w^*$ is a global minimizer of $\mathcal{E}$ (and thus a weak solution to \eqref{eq:main_problem-1dim}).

    Finally, if \eqref{ass:f>=} is replaced with \eqref{ass:f>}, then $F(x,\cdot)$ is strictly increasing on $[0,\infty)$ for a.e. $x\in (-1,1)$. Then, for any global minimizer $u \in W^{1,p}_0(-1,1)$ of $\mathcal{E}$, taking $v = u|_{[-1,0]}$ and using \eqref{eq:dichotomy} we get either $|v| = \overline{v}$ or
    \begin{equation}
        \label{eq:inegralfromx}
        \int_{x^*}^0 F\big(\xi, \overline{v}(\xi)\big)\,d\xi > \int_{x^*}^0 F\big(\xi, |v(\xi)|\big)\,d\xi = \int_{x^*}^0 F\big(\xi, v(\xi)\big)\,d\xi.
    \end{equation}
    However, \eqref{eq:inegralfromx} combined with \eqref{eq:J_Jhat_chain} immediately gives $\mathcal{J}(\overline{v}) = \widehat{\mathcal{J}}(\overline{v}) < \widehat{\mathcal{J}}(v) \leq \mathcal{J}(v)$, which contradicts \eqref{eq:u1_minimizes_J}. Consequently, $|v(x)| = \overline{v}(x) \geq 1 + x > 0$ for every $x \in (-1,0]$ and $|v'(x)| = \overline{v}'(x) \geq 1$ for a.e. $x \in (-1,0)$. Applying the same argument to the function $v_2(x) = u(-x)$ for $x \in [-1,0]$, which also minimizes $\mathcal{J}$ by Lemma~\ref{lem:comparison_of_minimizers}, we obtain the corresponding inequalities on $[0,1)$. Therefore, $|u(x)| \geq 1 - |x| > 0$ for every $x \in (-1,1)$ and $|u'(x)| \geq 1$ for a.e. $x \in (-1,1)$. Since $u(-1) = u(1) = 0$, Rolle's theorem implies $u \notin C^1(-1,1)$.
\end{proof}
We are finally in a position to state a corollary on the weak solvability of the one-dimensional analogue of the problem considered in Remark~\ref{remark1}.
\begin{remark}
    Let $p > 1$ and $q,r \in (1,p)$. The problem
    \begin{equation}
        \label{eq:problem1D-pqr}
        \begin{cases}
            -\left(|u'(x)|^{p-2}u'(x) - |u'(x)|^{q-2}u'(x)\right)' = |u(x)|^{r-2}u(x) & \text{for } x \in (-1,1), \\
            u(-1) = u(1) = 0 & 
        \end{cases}
    \end{equation}
    has at least one nonzero and nonnegative weak solution $u^*\in W^{1,p}_0(-1,1) \setminus C^1(-1,1)$ satisfying 
    \begin{equation*}
        \begin{split}
            \tfrac{1}{p}\int_{-1}^1 \left|\left(u^*\right)'(x)\right|^p dx & - \tfrac{1}{q}\int_{-1}^1 \left|\left(u^*\right)'(x)\right|^q dx - \tfrac{1}{r}\int_{-1}^1 \left|u^*(x)\right|^r dx \\
            & = \min_{u\in W^{1,p}_0(-1,1)} \left(\tfrac{1}{p}\int_{-1}^1 \left|u'(x)\right|^p dx - \tfrac{1}{q}\int_{-1}^1 \left|u'(x)\right|^q dx - \tfrac{1}{r}\int_{-1}^1 \left|u(x)\right|^r dx\right).
        \end{split}
    \end{equation*}
\end{remark}
\subsection{Final remark on multidimensional extensions}
In the reasoning above, the one-dimensional nature of the problem played a crucial role, and it is not clear whether the proposed approach can be directly extended to the general multidimensional setting. Nevertheless, the method employed here may be viewed as a one-dimensional instance of a more general technique available in the literature, which applies to multidimensional counterparts of equation~\eqref{eq:problem1D-pqr} and allows for considerably more general functions $\Phi$, provided that they do not depend explicitly on~the spatial variable. In particular, the aforementioned \cite[Theorem~4.2]{Kromer} implies that the problem
\begin{equation*}
    \begin{cases}
        -\Delta_p u(\boldsymbol{x}) + \Delta_q u(\boldsymbol{x}) = |u(\boldsymbol{x})|^{r-2}u(\boldsymbol{x}) & \text{for } \boldsymbol{x} \in B_1, \\
        u(\boldsymbol{x}) = 0 & \text{for } \boldsymbol{x} \in \partial B_1,
    \end{cases}
\end{equation*}
where $B_1 = \{\boldsymbol{x} \in \mathbb{R}^d : |\boldsymbol{x}| < 1\}$, admits a solution that is a global minimizer of the energy functional $\mathcal{E}\colon W^{1,p}_0(B_1) \longrightarrow \mathbb{R}$ given by
\begin{equation*}
    \mathcal{E}(u) = \tfrac{1}{p}\int_{B_1} |\nabla u(\boldsymbol{x})|^p \, d\boldsymbol{x} - \tfrac{1}{q}\int_{B_1} |\nabla u(\boldsymbol{x})|^q \, d\boldsymbol{x} - \tfrac{1}{r}\int_{B_1} |u(\boldsymbol{x})|^r \, d\boldsymbol{x}\quad \text{for every }u\in W^{1,p}_0(B_1).
\end{equation*}
An argument analogous to that used in the proof of Theorem~\ref{thm:generalized_existence} shows that this solution is nontrivial. Moreover, \cite[Theorem~4.2]{Kromer} ensures the existence of a radially symmetric minimizer of $\mathcal{E}$ that does not change sign in $B_1$, whereas \cite[Corollary~4.4]{Kromer} shows that such a minimizer is not of class $C^1$. It~remains unclear, however, whether this more general technique can be adapted to nonlinearities $f$ that depend explicitly on $\boldsymbol{x}$.




\end{document}